\documentclass[12pt]{amsart}
\usepackage{amsfonts,color,latexsym,amssymb,amsmath,epsfig,rotating,array}

\newtheorem{kapra}{Kapranov's Non-Archimedean Amoeba Theorem (special case)}

\newtheorem{adelic}{Adelic SPS-Conjecture}

\newtheorem{dfn}{Definition}[section]
\newtheorem{rem}[dfn]{Remark} 
\newtheorem{prop}[dfn]{Proposition}

\newtheorem{thm}[dfn]{Theorem} 
 
\newtheorem{lemma}[dfn]{Lemma}
\newtheorem{cor}[dfn]{Corollary}

\newtheorem{ex}[dfn]{Example}
\definecolor{purple}{rgb}{.5,0,.5}
\definecolor{red}{rgb}{.6,0,0} 
\definecolor{green}{rgb}{0,.5,0}

\renewcommand{\qed}{$\blacksquare$}

\newcommand{\sps}{\mathrm{SPS}}
\newcommand{\floor}[1]{\left\lfloor#1\right\rfloor}

\newlength{\burg}
\newlength{\koi}
\newlength{\sma}
\newlength{\jmr}
\newcommand{\thth}{^{\text{\underline{th}}}} 
\newcommand{\nd}{^{\text{\underline{nd}}}}

\newcommand{\stst}{^{\text{\underline{st}}}}
\newcommand{\ord}{\mathrm{ord}}

\newcommand{\np}{{\mathbf{NP}}}
\newcommand{\ch}{{\mathbf{CH}}}
\newcommand{\sharpp}{{\mathbf{\# P}}}

\newcommand{\glm}{\mathbb{G}\mathbb{L}_m}

\newcommand{\pp}{\mathbf{P}}
\newcommand{\vp}{\mathbf{VP}}
\newcommand{\vnp}{\mathbf{VNP}}
\newcommand{\ph}{\mathbf{PH}}

\newcommand{\pspa}{\mathbf{PSPACE}}

\newcommand{\eps}{\varepsilon}

\newcommand{\Q}{\mathbb{Q}}
\newcommand{\R}{\mathbb{R}}

\newcommand{\C}{\mathbb{C}}

\newcommand{\N}{\mathbb{N}}
\newcommand{\Z}{\mathbb{Z}}
\newcommand{\bO}{\mathbf{O}}

\newcommand{\Zn}{\Z^n}

\newcommand{\Qn}{\Q^n}
\newcommand{\Rn}{\R^n}

\newcommand{\Cn}{\C^n}

\newcommand{\Cs}{\C^*}
\newcommand{\Ks}{K^*}

\newcommand{\cR}{{\mathcal{R}}}
\newcommand{\cV}{{\mathcal{V}}}
\newcommand{\blah}{{\overline{\cV}}}
\newcommand{\klah}{{\overline{\cR}}}

\newcommand{\dia}{$\diamond$}
\newcommand{\trop}{\mathrm{Trop}}

\newcommand{\newt}{\mathrm{Newt}}

\newcommand{\supp}{\mathrm{Supp}}

\author{Pascal Koiran}
\email{pascal.koiran@ens-lyon.fr} 
\author{Natacha Portier} 
\email{natacha.portier@gmail.com} 
\author{J.\ Maurice Rojas}
\email{rojas@math.tamu.edu} 
\thanks{P.K.\ and N.P.\ were partially supported by the 
European Community (7th PCRD Contract: PIOF-GA-2009-236197). 
J.M.R.\ was partially supported by NSF MCS grant DMS-0915245 
and Labex MILYON}  
\title[Degenerate Valuations and the Hardness of the Permanent]{\mbox{}\\  
\vspace{-.25in}Counting Tropically Degenerate Valuations and 
\scalebox{1}[1]{$p$}-adic 
Approaches to the Hardness of the Permanent}  

\begin{document}

\mbox{} 
\vspace{-.5in} 
\begin{abstract}  
The Shub-Smale $\tau$-Conjecture is a hitherto unproven statement (on 
integer roots of polynomials) whose truth implies both a variant of 
$\pp\!\neq\!\np$ (for the BSS model over $\C$) and the hardness of 
the permanent. We give alternative conjectures, some 
potentially easier to prove, whose truth still implies the hardness of the 
permanent. Along the way, we discuss new upper bounds on the number of 
$p$-adic valuations of roots of certain sparse polynomial systems, 
culminating in a connection between quantitative $p$-adic geometry and 
complexity theory.
\end{abstract} 

\maketitle 

\vspace{-.4cm}
\mbox{}\hfill {\bf Dedicated to Mike Shub, on his 70$\thth$ birthday.} 
\hfill \mbox{}\\

\vspace{-.4cm}
\section{Introduction}  
Deep questions from algebraic complexity, cryptology, and 
arithmetic geometry can be approached 
through sufficiently sharp upper bounds on the number of roots of structured 
polynomials in one variable. (We review four such results 
in Section \ref{sub:earlier} below.) The main focus of this paper is the 
connection between the number of distinct {\em norms} of roots of polynomials, 
over the $p$-adic rationals $\Q_p$, and separations of complexity classes. Our 
first main theorem motivates the introduction of $p$-adic methods.  
\begin{dfn}
\label{dfn:sps} 
We define $\sps(k,m,t)$ 
to be the family of polynomials presented in 
the form $\sum^k_{i=1} \prod^m_{j=1} f_{i,j}$ where, for all $i$ and $j$, 
$f_{i,j}\!\in\!\Z[x_1]\!\setminus\!\{0\}$ and has at most 
$t$ monomial terms. We call such polynomials {\em SPS} (for {\em 
sum-product-sparse}) polynomials. \dia
\end{dfn}
\begin{thm}
\label{thm:koi}
Suppose that there is a prime $p$ with the following property:  
For all $k,m,t\!\in\!\N$ and $f\!\in\!\sps(k,m,t)$, 
we have that the cardinality of\\
\mbox{}\hfill $\displaystyle{S_f:=\left\{e\!\in\!\N \; : \; 
x\!\in\!\Z,\; f(x)\!=\!0,\; p^e|x, \text{ and } 
p^{e+1}\!\!\not|x. \right\}}$
\hfill\mbox{}\\ 
is $(kmt)^{O(1)}$. Then the permanent of $n\times n$ matrices cannot be 
computed by constant-free, division-free arithmetic circuits of size 
$n^{O(1)}$. 
\end{thm}
\begin{rem} 
The hypothesis of Theorem \ref{thm:koi} was the inspiration for  
this paper, since it is easily implied by the famous 
Shub-Smale $\tau$-Conjecture (see \cite{21a,21b} and Section 
\ref{sub:earlier} below).  
Theorem \ref{thm:koi} can in fact be strengthened further by weakening its 
hypothesis in various ways: see Remarks \ref{rem:stronger} and \ref{rem:best} 
of Section \ref{sec:koi} below. \dia 
\end{rem}  

\noindent 
The special cases of the hypothesis where $k=1$, $t=1$, or $m$ is a fixed 
constant are easy to prove (see, e.g., Lemma \ref{lemma:newt} of 
Section \ref{sub:uniup} below). 
However,  the hypothesis already becomes an open problem for $k=2$ or $t=2$.
The greatest $e$ such that $p^e$ divides an integer $x$ is nothing 
more than the {\em $p$-adic valuation} of $x$, hence our focus on $p$-adic 
techniques. 

We now describe certain families of univariate polynomials, 
and multivariate polynomial systems, where valuation counts in the 
direction of Theorem \ref{thm:koi} can actually be proved. In particular, we 
give another related hypothesis (in Theorem 
\ref{thm:hard} below), entirely within the realm of $p$-adic geometry, 
whose truth also implies the hardness of the permanent. 

\subsection{Provable Upper Bounds on the Number of Valuations} 
\begin{dfn} 
For any commutative ring $R$, we let $R^*\!:=\!R\!\setminus\!\{0\}$. The 
{\em support} of a polynomial $f\!\in\!R\!\left[x^{\pm 1}_1,\ldots,x^{\pm 1}_n
\right]$, denoted $\supp(f)$, is the set of exponent vectors appearing in the 
monomial term expansion of $f$. For any prime $p$ and $x\!\in\!\Z\setminus
\{0\}$ we let $\ord_p(x)$ denote the $p$-adic valuation of $x$, and we set 
$\ord_p(0)\!:=\!+\infty$. We then set $\ord_p(x/y)\!:=\!
\ord_p(x)-\ord_p(y)$ to extend $\ord_p(\cdot)$ to $\Q$, and we let $\Q_p$ 
denote the completion of $\Q$ with 
\scalebox{1}[1]{respect to the metric defined by
$|u-v|_p\!:=\!p^{-\ord_p(u-v)}$. The {\em $p$-adic complex numbers}, $\C_p$, 
are}\linebreak 
then the elements of the completion of the algebraic closure of $\Q_p$. 
Finally, for any polynomials $f_1,\ldots,f_r\!\in\!R\!\left[x^{\pm 1}_1,\ldots,
x^{\pm 1}_n\right]$, we let $Z_R(f_1,\ldots, f_r)$ (resp.\ 
$Z^*_R(f_1,\ldots,f_r)$) denote the set of 
\scalebox{.96}[1]{roots of $(f_1,\ldots,f_r)$ in $R^n$ (resp.\ 
$(R^*)^n$), and we use $\#S$ to denote the cardinality of a set $S$. \dia}  
\end{dfn} 

\noindent 
In particular, $\ord_p(\cdot)$ and $|\cdot|_p$ extend naturally to 
$\C_p$, and the algebraic closure of $\Q$ embeds naturally within $\C_p$. 
\cite{artin,weiss,serre,schikhof,robert,gouvea,katok} are some excellent 
sources for further background on $p$-adic fields. What will be most important 
for our setting is that $p$-adic norms (or, equivalently, $p$-adic valuations)
enable new hypotheses --- closer to being provable with current  
techniques --- that imply new separations of complexity classes. 

We will ultimately focus on counting valuations of roots of 
polynomial systems with few monomial terms as a means of 
understanding the valuations of roots of univariate SPS polynomials. 
For example, a simple consequence of our main multivariate 
bounds (Theorems \ref{thm:upper} and \ref{thm:mult} below) is the following 
univariate bound revealing that at least part of the $k\!=\!t\!=\!2$ case of 
the hypothesis of Theorem \ref{thm:koi} is true. 
\begin{cor} 
\label{cor:upper} 
Suppose $m_1,m_2\!\in\!\N$; $\alpha_i,\beta_i\!\in\!\C_p$; 
$\gamma_{i,j}\!\in\!\Z$; \\
\mbox{}\hfill $\displaystyle{f(x_1)\!:=\!\left(\prod^{m_1}_{i=1}(\alpha_{i,1}+
\beta_{i,1}x_1)^{\gamma_{i,1}}\right)+\left(\prod^{m_2}_{i=1}(\alpha_{i,2}+
\beta_{i,2}x_1)^{\gamma_{i,2}}\right)}$\hfill\mbox{}\\ 
is not identically zero; and the lower hulls of the $p$-adic Newton polygons 
(cf.\ Definition \ref{dfn:newt} below) of the two products have 
no common vertices. 
Then $\#\ord_p\!\left(Z_{\C_p}\!\left(f\right)\right)\!\leq\!m_1+m_2$, 
and this bound is tight. Furthermore, any root of $f$ in $\C_p$ not 
making both products vanish has multiplicity at most $m_1+m_2$, and this 
bound is tight as well. \qed 
\end{cor}  

\noindent 
Bounds for the number of valuations, independent of the degree, 
had previously been known only for sparse polynomials, i.e., 
polynomials in $\sps(k,1,t)$: 
see, e.g., Lemma \ref{lemma:newt} of Section \ref{sub:uniup} and  
\cite{weiss}. 
Note in particular that $\sps(2,m,2)$ contains the family of $f$ in our 
corollary when $\gamma_{i,j}\!=\!1$ for 
all $i,j$. Also, our valuation count above is independent of the 
$\alpha_i, \beta_i,\gamma_i$.  

\noindent 
\vbox{\scalebox{.93}[1]{We say that $F$ is {\em tropically generic (over 
$\C_p$)} iff the closures of $\ord_p(Z^*_{\C_p}(f_1)),\ldots,
\ord_p(Z^*_{\C_p}(f_n))$} 
intersect transversally.  {\em Kapranov's Non-Archimedean Theorem} 
\cite{kapranov}, reviewed in Section \ref{sec:back} below, tells us that the 
closure of each $\ord_p(Z^*_{\C_p}(f_i))$ is in fact a polyhedral complex of 
codimension $1$ in $\Rn$, so it makes sense to speak of transversality.  

\noindent 
\begin{picture}(200,50)(5,-88)
\put(370,-105){
\epsfig{file=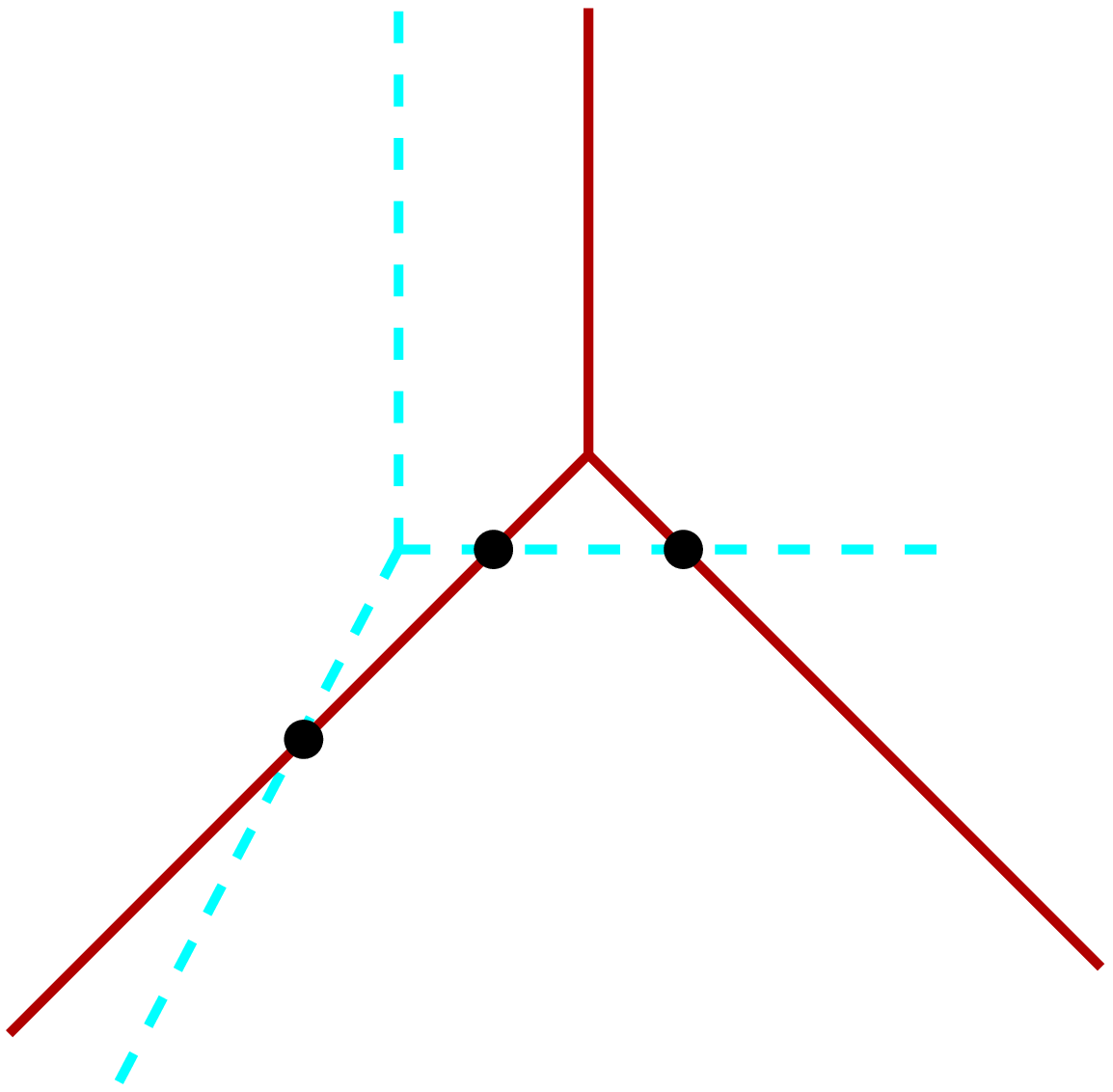,height=1.1in,clip=}
}
\put(0,-95){
\begin{minipage}[b]{1\linewidth}
\vspace{0pt}
\begin{ex} 
\scalebox{.95}[1]{For any prime $p$, the polynomials 
$f_1\!:=\!x_1x_2-p-x^2_1$ and}\\ 
\scalebox{1}[1]{$f_2\!:=\!x_2-1-px^2_1$ have 
$\ord_p(Z^*_{\C_p}(f_1))$ and $\ord_p(Z^*_{\C_p}(f_2))$ intersecting}\\  
\scalebox{.93}[1]{transversally as shown on the right. 
$\ord_p(Z^*_{\C_p}(f_1))$ (resp.\ $\ord_p(Z^*_{\C_p}(f_2))$)}\\  
\scalebox{.96}[1]{consists of the rational points on the 
solid (resp.\ dashed) curve. \dia} 
\end{ex} 
\end{minipage} 
} 
\end{picture} 
}

\vspace{.5cm}
\noindent
\begin{picture}(200,65)(5,-100)
\put(380,-102){
\epsfig{file=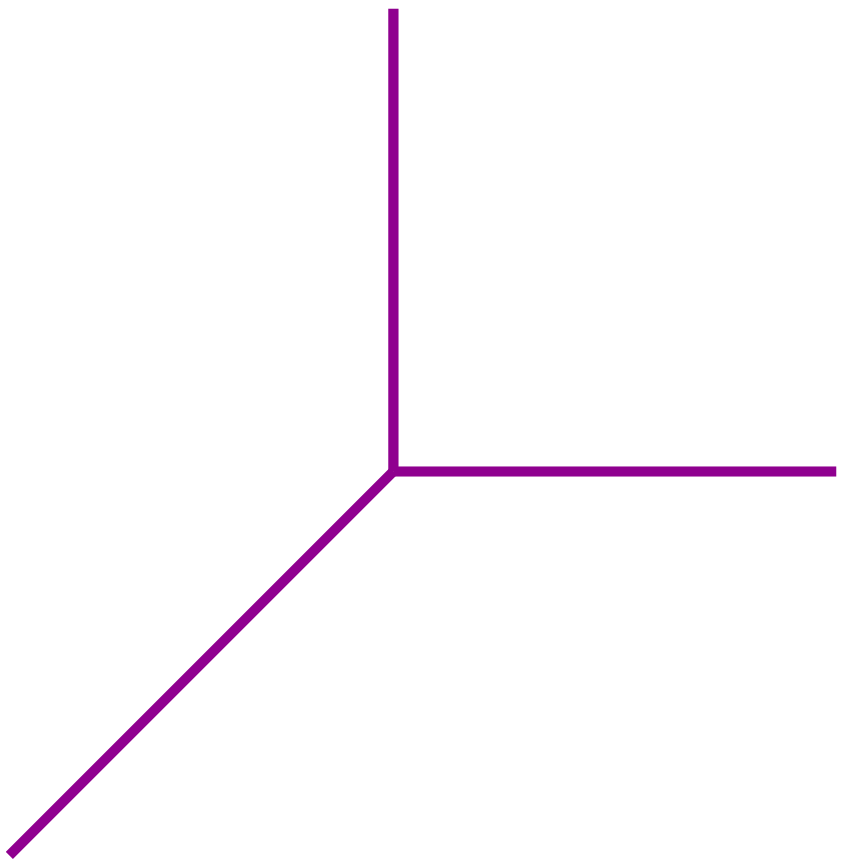,height=1.1in,clip=}
}
\put(0,-100){
\begin{minipage}[b]{1\linewidth}
\vspace{0pt}
\begin{ex}
For any prime $p$, the system $F\!:=\!(x+y+1,x+y+1+p)$\\ 
\scalebox{.97}[1]{shows us that having just finitely many roots over $\C_p$ 
need {\em not} imply tropical}\\ 
genericity. In particular, while $F$ has no roots at all in $\C^2_p$, we 
have that\\ 
$\ord_p(Z^*_{\C_p}(x+y+1)), \ord_p(Z^*_{\C_p}(x+y+1+p))\!\subset\!\R^2$ are 
identical and\\ 
exactly the set of rational points on the right-hand union of $3$ rays:  
\end{ex} 
\end{minipage} 
}
\end{picture} 

\noindent 
{\em (So the intersection is non-transversal.) 
Nevertheless, $\ord_p(Z^*_{\C_p}(F))$ is empty. \dia }  

\vspace{-.1cm} 
\begin{rem} 
Throughout our paper, we abuse notation slightly by 
setting $\ord_p(x_1,\ldots,x_n)\!:=\!(\ord_p(x_1),\ldots,\ord_p(x_n))$. 
The set $\ord_p\!\left(Z^*_{\C_p}(f_1,\ldots,f_r)\right)$ is thus well-defined 
but, as revealed above, need {\em not} be the same as 
$\bigcap^r_{i=1} \ord_p\!\left(Z^*_{\C_p}(f_i)\right)$. \dia 
\end{rem} 
\begin{dfn}
For any finite subsets $A_1,\ldots,A_n\!\subset\!\Zn$ we define 
$\blah_p(A_1,\ldots,A_n)$ (resp.\ $\klah_p(A_,\ldots,A_n)$) to be the maximum 
of $\#\ord_p\!\left(Z^*_{\C_p}(F)\right)$ (resp.\ 
$\ord_p\!\left(Z^*_{\Q_p}(F)\right)$) over 
all $F\!:=$\linebreak
\scalebox{.94}[1]{$(f_1,\ldots,f_n)$ with $f_i\!\in\!\C_p[x_1,\ldots,x_n]$ and  
$\supp(f_i)\!\subseteq\!A_i$ for all $i$, and $Z^*_{\C_p}(F)$ finite.
We also define}\linebreak 
\scalebox{.97}[1]{$\cV_p$ to be the corresponding analogue of 
$\blah_p$ where we restrict further to tropically generic $F$. \dia}  
\end{dfn} 

\noindent 
\scalebox{.95}[1]{Clearly 
$\cV_p(A_1,\ldots,A_n)\!\leq\!\blah_p(A_1,\ldots,A_n)$ and 
$\klah_p(A_1,\ldots,A_n)\!\leq\!\blah_p(A_1,\ldots,A_n)$. 
While {\em Smirnov's}}\linebreak 
{\em Theorem} \cite[Thm.\ 3.4]{smirnov} 
implies that $\cV_p(A_1,\ldots,A_n)$ is well-defined and finite for any fixed 
$(A_1,\ldots,A_n)$, explicit upper bounds for $\blah_p(A_1,\ldots,A_n)$ appear 
to be unknown. So we derive such an upper bound for certain $(A_1,\ldots,A_n)$. 
\begin{thm}
\label{thm:upper} 
Suppose $p$ is any prime, $A_1,\ldots,A_n\!\subset\!\Zn$, 
$A\!:=\#\bigcup_i A_i$, $t\!:=\!\#A$, and $e_i$ denotes the $i\thth$ standard 
basis vector of $\Rn$. Then:  
\begin{enumerate}
\addtocounter{enumi}{-1} 
\item{\ \scalebox{.86}[1]{
$t\leq n \ \ \ \ \; \Longrightarrow 
\cV_p(A_1,\ldots,A_n)=\blah_p(A_1,\ldots,A_n)=0$}. } 
\item{\ 
\scalebox{.86}[1]{
$t=n+1 \Longrightarrow 
\cV_p(A_1,\ldots,A_n)=\blah_p(A_1,\ldots,A_n)\leq 1$. In 
particular, $\cV_p\!\left(\{\bO,e_1\},\ldots,\{\bO,e_n \}\right)\!=\!1$.}} 
\item{\; \scalebox{.97}[1]{[$t=n+2$ and every collection of $n$ distinct pairs 
of points of $A$ determines an $n$-tuple}\linebreak  
of linearly independent vectors] $\Longrightarrow 
\blah_p(A_1,\ldots,A_n)\leq \max\left\{2,\floor{\frac{n}{2}}^n+n\right\}$. 
Also, 
$\cV_p\!\left(\{\bO,2e_1,e_1+e_2\},\{\bO,2e_1,e_2+e_3\},
\ldots,\{\bO,2e_1,e_{n-1}+e_n\},\{\bO,2e_1,e_n\}\right)\!=\!n+1$.}   
\end{enumerate} 
\end{thm} 

\noindent 
We conjecture that the upper bound in Assertion (2) can in fact 
be improved to $n+1$. 
It is easily shown that the general position assumption on $A$ holds for a 
dense open set of exponents. For instance, if $A$ has convex 
hull an $n$-simplex then the hypothesis of Assertion (2) holds automatically. 

Whether the equality $\cV_p(A_1,\ldots,A_n)=\blah_p(A_1,\ldots,A_n)$ 
holds beyond the setting of Assertions (0) and (1) is an intriguing open 
question. In fact, proving just that the 
growth of orders of $\klah_p(A_1,\ldots,A_n)$ and 
$\cV_p(A_1,\ldots,A_n)$ differ by a constant has deep implications.  
\begin{thm} 
\label{thm:hard} 
Suppose there is a prime $p$ such that 
$\klah_p(A_1,\ldots,A_n)\!=\!\cV_p(A_1,\ldots,A_n)^{O(1)}$
for all finite $A_1,\ldots,A_n\!\subset\!\Zn$. Then the permanent of 
$n\times n$ matrices cannot be computed by constant-free, division-free 
arithmetic circuits of size $n^{O(1)}$. 
\end{thm} 

\noindent 
We thus obtain an entirely tropical geometric statement implying the 
hardness of the\linebreak 
permanent. In Theorem \ref{thm:hard} it in fact suffices to restrict to 
certain families of supports $A_i$ (see Proposition \ref{prop:red} below). 
We are currently unaware of any $(A_1,\ldots,A_n)$ {\em not} satisfying the 
equality $\cV_p(A_1,\ldots,A_n)\!=\!\blah_p(A_1,\ldots,A_n)$ for all primes 
$p$. 

Intersection multiplicity is a key subtlety underlying the counting of 
valuations. 
\begin{thm} 
\label{thm:mult} 
Suppose $K$ is any algebraically closed field of characteristic $0$ 
and $A\!\subset\!\Zn$ has cardinality at most $n+2$ and no $n+1$ points 
of $A$ lie in a hyperplane. 
Suppose also that $f_1,\ldots,
f_n\!\in\!K[x_1,\ldots,x_n]$ have support contained in $A$ and 
$\#Z^*_{K}(F)$ is finite. Then the intersection multiplicity of any 
point of $Z^*_{K}(F)$ is at most $n+1$ (resp.\ $1$) when $\#A\!=\!n+2$ 
(resp.\ $\#A\!=\!n+1$), and both bounds are sharp.  
\end{thm} 

\noindent 
The intersection multiplicity considered in our last theorem is the 
classical definition coming from commutative algebra or differential 
topology (see, e.g., \cite{ifulton}). 

\subsection{Earlier Applications of Root Counts for Univariate 
Structured Polynomials} 
\label{sub:earlier} 
Recall the following classical definitions on the evaluation complexity 
of univariate polynomials. 
\begin{dfn}
For any field $K$ and $f\!\in\!K[x_1]$ let $s(f)$ --- the {\em SLP
complexity of $f$} --- denote the smallest $n$
such that $f\!=\!f_n$ identically where the sequence $(f_{-N},\ldots,f_{-1},
f_0,\ldots,f_n)$ satisfies the following conditions:
$f_{-1},\ldots,f_{-N}\!\in\!K$, $f_0\!:=\!x_1$, and, for all $i\!\geq\!1$,
$f_i$ is a sum, difference, or product of some pair of elements
$(f_j,f_k)$ with $j,k\!<\!i$. Finally, for any $f\!\in\!\Z[x_1]$, we
let $\tau(f)$ denote the obvious analogue of $s(f)$ where the definition
is further restricted by assuming $N\!=\!1$ and $f_{-1}\!:=\!1$. \dia
\end{dfn}

\noindent
Note that we always have $s(f)\!\leq\!\tau(f)$ since $s$ does not
count the cost of computing large integers (or any constants). 
One in fact has $\tau(n)\!\leq\!2\log_2 n$ for any
$n\!\in\!\N$ \cite[Prop.\ 1]{svaiter}.
See also \cite{brauer,moreira} for further background. 

We can then summarize some seminal results of 
B\"urgisser, Cheng, Lipton, Shub, and Smale as follows:  
\begin{thm}
\label{thm:tau} \mbox{}\\
I. (See \cite[Thm.\ 3, Pg.\ 127]{bcss} and \cite[Thm.\ 1.1]{burgtau}.)
Suppose that for all nonzero $f$\linebreak
\mbox{}\hspace{.5cm}\scalebox{.94}[1]{$\in\!\Z[x_1]$ we have
$\#Z_\Z(f)\!\leq\!\tau(f)^{O(1)}$.
Then (a) $\pp_\C\!\neq\!\np_\C$ and (b) the permanent of $n\times n$ 
}\linebreak
\mbox{}\hspace{.5cm}\scalebox{.97}[1]{matrices cannot be computed by 
constant-free, division-free arithmetic circuits of size $n^{O(1)}$.} 

\noindent
II. (Weak inverse to (I) \cite{lipton}.\footnote{
Lipton's main result from \cite{lipton} is in fact stronger, allowing for
rational roots and primes with a mildly differing number of digits.}) If
there is an $\eps\!>\!0$ and a sequence $(f_n)_{n\in\N}$
of polynomials\linebreak
\mbox{}\hspace{.7cm}in $\Z[x_1]$ satisfying:\\
\mbox{}\hspace{1.2cm}(a) $\#Z_\Z(f_n)\!>\!e^{\tau(f_n)^\eps}$ for all
$n\!\geq\!1$ \ and  \
(b) $\deg f_n, \max\limits_{\zeta\in Z_\Z(f)}|\zeta|\!\leq\!
2^{(\log \#Z_\Z(f_n))^{O(1)}}$\\
\mbox{}\hspace{.7cm}then, for infinitely many $n$, at least
$\frac{1}{n^{O(1)}}$ of the $n$ digit integers that are products of
exactly\linebreak
\mbox{}\hspace{.7cm}two distinct primes (with an equal number of
digits) can be factored by a Boolean circuit\linebreak
\mbox{}\hspace{.7cm}of size $n^{O(1)}$.

\noindent
III. (Number field analogue of (I) implies Uniform Boundedness \cite{cheng}.)
Suppose that for\linebreak
\mbox{}\hspace{.8cm}any number field $K$ and $f\!\in\!K[x_1]$ we have
$\#Z_K(f)\!\leq\!c_1 1.0096^{s(f)}$, with $c_1$ depending only\linebreak
\mbox{}\hspace{.8cm}on $[K:\Q]$. Then there is a constant $c_2\!\in\!\N$
depending only on $[K:\Q]$ such that for any\linebreak
\mbox{}\hspace{.8cm}elliptic curve $E$ over $K$, the torsion
subgroup of $E(K)$ has order at most $c_2$. \qed
\end{thm}

\noindent
The hypothesis in Part (I) is known as the {\em (Shub-Smale)
$\tau$-Conjecture} and was stated as the fourth 
problem (still unsolved as of late 2013) on Smale's list of the most important 
problems for the $21\stst$ century \cite{21a,21b}.
Via fast multipoint evaluation applied to the polynomial 
$(x-1)\cdots(x-m^2)$ \cite{compualg} one can show that the 
$O$-constant from the $\tau$-Conjecture should be at least
$2$ if the $\tau$-Conjecture is true. 

The complexity classes $\pp_\C$ and $\np_\C$ are respective analogues
(for the BSS model over $\C$ \cite{bcss}) of the well-known complexity
classes $\pp$ and $\np$.  (Just as in the famous $\pp$ vs.\ $\np$ Problem, the 
equality of $\pp_\C$ and $\np_\C$ remains an open question.)
The assertion on the hardness of the permanent in Theorem \ref{thm:tau}
is also an open problem and its proof would be a major step toward
solving the {\em $\vp$ vs.\ $\vnp$ Problem} ---
Valiant's algebraic circuit analogue of the $\pp$ vs.\ $\np$ 
Problem~\cite{valiant,burgcook,koiran,jml}: 
The only remaining issue to resolve for a complete solution of this problem 
would then be the restriction to constant-free circuits in Part (I).  

The hypothesis of Part (II) (also unproved as of late 2013) merely posits a 
sequence of polynomials
violating the $\tau$-Conjecture in a weakly exponential manner.
The conclusion in Part (II) would violate a widely-believed
version of the cryptographic hardness of integer factorization.

The conclusion in Part (III) is the famous {\em Uniform Boundedness
Theorem}, due to Merel \cite{merel}.
Cheng's conditional proof (see
\cite[Sec.\ 5]{cheng}) is dramatically simpler and would yield effective bounds
significantly improving known results (e.g., those of Parent \cite{parent}).
In particular, the $K\!=\!\Q$ case of the hypothesis of Part (III)
would yield a new proof (less than a page
long) of Mazur's landmark result on torsion points \cite{mazur}.

More recently, Koiran has suggested real analytic methods (i.e., 
upper bounds on the number of real roots)  
as a means of establishing the desired upper bounds on the number of integer 
roots \cite{koiran}, and Rojas has suggested $p$-adic methods \cite{adelic}. 
In particular, the following variation on the hypothesis from 
Theorem \ref{thm:koi} appears in slightly more refined form in \cite{adelic}: 
\begin{adelic}
For any $k,m,t\!\in\!\N$ and $f\!\in\!\sps(k,m,t)$,
there is a field\linebreak
\scalebox{1}[1]{$L\!\in\!\{\R,\Q_2,\Q_3,\Q_5,\dots\}$ such that
$f$ has no more than $(kmt)^{O(1)}$ distinct roots in $L$.}
\end{adelic}  
\begin{thm} 
\label{thm:adelic} 
If the Adelic SPS-Conjecture is true then the permanent of $n\times n$ 
matrices cannot be computed by constant-free, division-free arithmetic 
circuits of size $n^{O(1)}$. 
\end{thm} 

\noindent 
{\bf Proof of Theorem \ref{thm:adelic}:} The truth of the Adelic 
SPS-Conjecture clearly implies the following special case of the 
Shub-Smale $\tau$-Conjecture: the number of integer roots of any 
$f\!\in\!\sps(k,m,t)$ is $(kmt)^{O(1)}$. The latter statement 
in turn implies the hypothesis of Theorem \ref{thm:koi},  
so by the conclusion of Theorem \ref{thm:koi} we are done. \qed 

\medskip 
Note that the Adelic SPS-Conjecture can not be simplified to counting 
just the valuations: {\em any} fixed polynomial in $\Z[x_1]\setminus\{0\}$ 
will have exactly 
{\em one} $p$-adic valuation for its roots in $\C_p$ for sufficiently large 
$p$. (This follows easily from, e.g., Lemma \ref{lemma:newt} of the next 
section.) An alternative simplification (and stronger hypothesis) would be to 
ask for a {\em single} field $L\!\in\!\{\R,\Q_2,\Q_3,\Q_5,\ldots
\}$ where the number of roots in $L$ of any $f\!\in\!\sps(k,m,t)$ is 
$(kmt)^{O(1)}$. The latter simplification is an open problem, although 
it is now known that one can not ask for too much more: the 
stronger statement that the number of roots in $L$ of any $f\!\in\!\Z[x_1]$ 
is $\tau(f)^{O(1)}$ is known to be false. Counter-examples 
are already known over $\R$ (see, e.g., \cite{bocook}), and over $\Q_p$ 
for any prime $p$ \cite[Example 2.5 \& Sec.\ 4.5]{adelic}.  

The latter examples are much more recent, so for the convenience of the 
reader we summarize them here: Recall that the $p$-adic {\em integers}, $\Z_p$, 
are those elements of $\Q_p$ with nonnegative valuation. 
(So $\Z\!\subsetneqq\!\Z_p$ in particular.) 
\begin{ex}
\label{ex:slp}
Consider the recurrence $h_1\!:=\!x_1(1-x_1)$
and $h_{n+1}\!:=\!\left(p^{3^{n-1}}-h_{n}\right)h_n$ for all
$n\!\geq\!1$. Then $h_n$ has
degree $2^n$, exactly $2^n$ roots in $\Z_p$, 
and $\tau(h_n)=O(n)$. However, the only integer roots of 
$h_n$ are $\{0,1\}$ (see \cite[Sec.\ 4.5]{adelic}). Note also that $h_n$ has 
just $n$ distinct valuations for its roots in $\C_p$. The last fact 
follows easily from Lemma \ref{lemma:newt}, stated in the next section. \dia 
\end{ex}

\noindent 
Note, however, that it is far from obvious if the polynomial $h_n$ above is  
in $\sps(k,m,t)$ for some triple $(k,m,t)$ of functions growing polynomially 
in $n$. 

\subsection{From Univariate SPS to Multivariate Sparse} 
\label{sub:uniup} 
Perhaps the simplest reduction of root counts for univariate SPS polynomial 
to root counts for multivariate sparse polynomial systems is the following. 
\begin{prop} 
\label{prop:red} 
Suppose 
$f\!\in\!\sps(k,m,t)$ 
is written 
$\sum^k_{i=1} \prod^m_{j=1} f_{i,j}$ as in Definition \ref{dfn:sps}. 
Let $F\!:=\!(f_1,\ldots,f_{km+1})$ be the polynomial system defined by 
$f_{km+1}(x_1,\ldots,y_{i,j},\ldots)\!:=\! 
\sum^k_{i=1} \prod^m_{j=1} y_{i,j}$ and 
$f_{(i-1)m+j}(x_1,y_{i,j})\!:=\!y_{i,j} - f_{i,j}(x_1)$ for all 
$(i,j)\!\in\!\{1,\ldots,k\}\times \{1,\ldots,m\}$.  
(Note that $F$ involves exactly 
$km+1$ variables; $f_1,\ldots,f_{km}$ each have at most $t+1$ 
monomial terms; and $f_{km+1}$ involves exactly $k$ monomial terms.)  
Then $f$ not identically zero implies that $F$ has only finitely many roots 
in $\C_p$, and the $x_1$-coordinates of the roots of $F$ in $\C_p$ 
are exactly the roots of $f$ in $\C_p$. \qed 
\end{prop} 

Upper bounds for the number of valuations of the roots of multivariate sparse 
polynomials can then, in some cases, yield useful upper bounds for the number 
of valuations of the roots of univariate SPS polynomials. 
\begin{lemma} 
\label{lemma:maybetrivial} 
Following the notation of Proposition \ref{prop:red}, suppose 
$F$ is tropically generic. Then $\#\ord_p\!\left(Z^*_{\C_p}(F)\right)\!\leq\! 
k(k-1)(2km(t-1)+1)/2\!=\!O(k^3mt)$. 
\end{lemma}  

\noindent 
The crux of our paper is whether the bound above holds without 
tropical genericity. We prove Lemma \ref{lemma:maybetrivial} 
in Section \ref{sec:back} below. 

To prove upper bounds such as Lemma \ref{lemma:maybetrivial} (and 
Corollary \ref{cor:upper} and Theorem \ref{thm:upper}) we will need to use 
some polyhedral geometric tricks. So let us first review $p$-adic Newton 
polygons (see, e.g., \cite{weiss,gouvea}).
\begin{dfn} 
\label{dfn:newt} 
Given any prime $p$ and a polynomial
$f(x_1)\!:=\!\sum^t_{i=1}c_ix^{a_i}_1$
$\in\!\C_p[x_1]$, we define its {\em $p$-adic Newton polygon},
$\newt_p(f)$, to be the convex hull of\footnote{i.e., smallest convex
set containing...} the points $\{(a_i,\ord_p c_i)\; | \; 
i\!\in\!\{1,\ldots,t\}\}$. Also, a face of a polygon $Q\!\subset\!\R^2$ is
called {\em lower} if and only if it has an inner normal with positive last
coordinate, and the {\em lower hull} of $Q$ is simply the union of all its
lower edges. Finally, the polynomial associated to summing the terms of $f$
corresponding to points of the form $(a_i,\ord_p c_i)$ lying on some lower
face of $\newt_p(f)$ is called a {\em ($p$-adic) lower polynomial}. \dia
\end{dfn}
\begin{ex}
\label{ex:newt}
For $f(x_1):=36 -8868x_1 +29305x^2_1 -35310x^3_1 +18240x^4_1
-3646x^5_1+243x^6_1$,\linebreak

\vspace{-.8cm}
\noindent
\begin{minipage}[t]{.8\linewidth}
\vspace{0pt}
the polygon $\newt_3(f)$ has exactly $3$ lower edges and
can easily be verified to resemble the
illustration to the right. The polynomial $f$ thus has
exactly $2$ lower binomials, and $1$ lower trinomial over $\C_3$. \dia
\end{minipage}\hspace{.6cm}
\begin{minipage}[t]{1in}
\vspace{0pt}
\epsfig{file=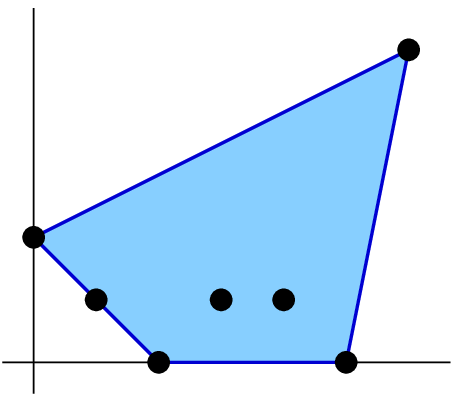,height=.6in}
\end{minipage} 
\end{ex}

The $p$-adic Newton polygon is particularly important because it 
allows us to count valuations (or norms) of $p$-adic complex roots exactly when 
the monomial term expansion is known. 
\begin{lemma}
\label{lemma:newt}
(See, e.g., \cite[Prop.\ 3.1.1]{weiss}.)
The number of roots of $f$ in $\C_p$ with valuation $v$, counting
multiplicities, is {\em exactly} the horizontal length of the lower face of
$\newt_p(f)$ with inner normal $(v,1)$. \qed
\end{lemma}
\begin{ex}
In Example \ref{ex:newt}, note that the $3$ lower
edges have respective horizontal lengths $2$, $3$, and $1$,
and inner normals $(1,1)$, $(0,1)$, and $(-5,1)$. Lemma \ref{lemma:newt}
then tells us that $f$ has exactly $6$ roots in $\C_3$:
$2$ with $3$-adic valuation $1$, $3$ with $3$-adic valuation $0$, and $1$
with $3$-adic valuation $-5$. Indeed, one can check that the roots
of $f$ are exactly $6$, $1$, and $\frac{1}{243}$, with respective
multiplicities $2$, $3$, and $1$. \dia
\end{ex}

To prove Theorem \ref{thm:hard} we will need to review 
the higher-dimensional version of the $p$-adic Newton polygon: 
the $p$-adic Newton {\em polytope}. 

\section{Background on $p$-adic Tropical Geometry} 
\label{sec:back}
The definitive extension of $p$-adic Newton polygons to arbitrary dimension 
(and general non-Archimedean, algebraically closed fields) is due to Kapranov. 
\begin{dfn}
\scalebox{.96}[1]{For any polynomial $f\!\in\!\C_p[x_1,\ldots,x_n]$ written 
$\sum_{a\in A} c_ax^a$ (with $x^a\!=\!x^{a_1}_1 \cdots x^{a_n}_n$}\linebreak  
understood) we define its {\em $p$-adic Newton polytope}, $\newt_p(f)$, 
to be the convex hull of the point set $\{(a,\ord_p(c_a))\; | \; 
a\!\in\!A\}$. We also define the {\em $p$-adic tropical variety of 
$f$} (or {\em $p$-adic amoeba of}\linebreak 
\scalebox{.905}[1]{$f$), $\trop_p(f)$, 
to be $\{v\!\in\!\Rn\; | \; (v,1) \text{ is an inner 
normal of a positive-dimensional face of } \newt_p(f)\}$. \dia}  
\end{dfn}

\noindent
We note that in \cite{kapranov}, the $p$-adic tropical variety of $f$ was 
defined via a {\em Legendre transform} (a.k.a.\ {\em support function} 
\cite{ziegler}) of the lower hull of $\newt_p(f)$. It is easy to see that both 
defintions are equivalent.
\begin{kapra}
\cite{kapranov}
Following the notation above, 
$\ord_p\!\left(Z^*_{\C_p}(f)\right)\!=\!\trop_p(f)\cap \Q^n$. \qed
\end{kapra}

A simple consequence of Kapranov's Theorem is that counting valuations is most
interesting for zero-dimensional algebraic sets. 
\begin{prop}
\label{prop:fin}
Suppose $f_1,\ldots,f_r\!\in\!\C_p\!\left[x^{\pm 1}_1,\ldots,x^{\pm 1}_n
\right]$, $F\!:=\!(f_1,\ldots,f_r)$, and $Z^*_{\C_p}(F)$ is infinite.
Then $\ord_p\!\left(Z^*_{\C_p}(F)\right)$ is infinite.
\end{prop}

\noindent
{\bf Proof:} By the definition of dimension for algebraic sets over
an algebraically closed field, there must be a 
linear projection $\pi : \Cn_p \longrightarrow I$, for some 
coordinate subspace $I$ of positive dimension $k$, with 
$\pi\!\left(Z^*_{\C_p}(F)\right)$ dense.
Taking valuations, and applying Kapranov's Theorem, 
this implies that $\ord_p\!\left(\pi\!\left(Z^*_{\C_p}(F)\right)\right)$
must be linearly isomorphic to $\Q^k$ minus a (codimension $1$) 
polyhedral complex. In other words, 
$\ord_p\!\left(Z^*_{\C_p}(F)\right)$ must be infinite. \qed 

Another consequence of Kapranov's Theorem is a simple characterization 
of $\ord_p(Z^*_{\C_p}(F))$ when $F\!:=\!(f_1,\ldots,f_n)$ is over-determined 
in a certain sense. This is based on a trick commonly used in 
toric geometry, ultimately reducing to an old matrix factorization: For any 
matrix $M=[M_{i,j}]\!\in\!\Z^{n\times n}$ and $x\!\in\!(\Cs_p)^n$, we define
$x^M\!:=\!\left(x^{M_{1,1}}_1\cdots x^{M_{n,1}},\ldots,
x^{M_{1,n}}_1\cdots x^{M_{n,n}}\right)$. We then call the map $m_M : 
(\Cs_p)^n \longrightarrow (\Cs_p)^n$ defined by $m_M(x)\!:=\!x^M$
a {\em monomial change of variables}.
\begin{lemma}
\label{lemma:mono}
Given any finite set $A\!=\!\{a_1,\ldots,a_n\}\!\subset\!\Zn$ lying in a 
hyperplane in $\Rn$, there is a matrix $U\!\in\!\Z^{n\times n}$, with 
determinant $\pm 1$, satisfying the following conditions:
\begin{enumerate}
\item{$Ua_i\!\in\!\Z^{i}\times\{0\}^{n-i}$ for all $i\!\in\!\{1,\ldots,n\}$.}
\item{Left (or right) multiplication by $U$ induces a linear bijection
of $\Zn$.}
\item{$m_U$ is an automorphism of the multiplicative group
$(\Cs_p)^n$, with inverse $m_{U^{-1}}$. In\linebreak 
particular, the map sending $\ord_p(x)\mapsto 
\ord_p(m_U(x))$ for all $x\!\in\!(\Cs_p)^n$ is a linear automorphism of 
$\Qn$. \qed}
\end{enumerate}
\end{lemma}

\noindent
Lemma \ref{lemma:mono} follows immediately from the existence
of {\em Hermite factorization} for matrices with integer entries
(see, e.g., \cite{hermite,storjophd}). In fact, the
matrix $U$ above can be constructed efficiently, but this need not
concern us here. The characterization of $\ord_p\!\left(Z^*_{\C_p}(F)\right)$ 
for over-determined $F$ is the following statement.  
\begin{prop} 
\label{prop:over} 
Suppose $A_1,\ldots,A_n\!\subseteq\!A\!\subset\!\Zn$ and $A$ lies 
in some $(n-1)$-flat of $\Rn$. Then $\cV_p(A_1,\ldots,A_n)\!=\!\blah_p
(A_1,\ldots,A_n)\!=\!0$. 
\end{prop} 

\noindent 
{\bf Proof:} Suppose $F\!:=\!(f_1,\ldots,f_n)$ where $\supp(f_i)\!\subseteq
\!A_i$ for all $i$. By Lemma \ref{lemma:mono} we may 
assume that $A\!\subset\!\Z^{n-k}\times \{0\}^k$ for some $k\!\geq\!1$.  
Clearly then, $Z^*_{\C_p}(F)$ is either empty or\linebreak 
contains a coordinate $k$-flat. So  
$\ord_p\!\left(Z^*_{\C_p}(F)\right)$ must either be empty or infinite,  
and we are done. \qed  

\medskip 
Another consequence of Kapranov's Theorem is the following characterization 
of\linebreak $\ord_p\!\left(Z^*_{\C_p}(f)\right)$ for certain trinomials. 
Recall that $\R_+$ is the set of positive real numbers and 
\scalebox{.97}[1]{that $\R_+v$, for any 
vector $v\!\in\!\R^N\setminus\{\bO\}$, is the {\em open ray} generated 
by all positive multiples of $v$.}  
\begin{lemma} 
\label{lemma:vert} 
Suppose $g\!\in\!\C_p[x_1]$ has exactly $t$ monomial terms, the 
lower hull of $\newt_p(g)$ consists of exactly $t'$ edges, and 
$f(x_1,y_i)\!:=\!y_i-g(x_1)$ is considered as a polynomial in 
$\C_p[x_1,y_1,\ldots,y_N]$ with $N\!\geq\!i$. Then $\ord_p(Z^*_{\C_p}(f))$ is 
the set of rational points of a polyhedral complex $\Sigma_f$ of the following 
form: a union of \textcolor{red}{{\bf (a)}} an open $(N-1)$-dimensional 
half-space parallel to $(\R_+(-1,\deg(g)))\times \R^{N-1}$, 
\textcolor{green}{{\bf (b)}} $t'$ ``vertical'' 
open half-spaces parallel to $(\R_+(0,1))\times \R^{N-1}$,
\textcolor{blue}{{\bf (c)}} $t'-1$ strips of the form $L\times \R^{N-2}$ where 
$L\!\subset\!\R^2$ is a line segment missing one of its vertices, and 
\textcolor{purple}{{\bf (d)}} a closed $(N-1)$-dimensional 
half-space parallel to $(\R_+\cup\{0\})\times \{0\} \times \R^{N-1}$. 
\end{lemma} 
\begin{ex} 
For any prime $p$, the polynomial \\ 
\mbox{}\hfill 
$f(x_1,y_1)\!:=\!y_1-(x^3_1-(1+p+p^2)x^2_1+(p+p^2+p^3)x_1
-p^3)$ \hfill \mbox{}\\ 
has $\ord_p\!\left(Z^*_{\C_p}(f)\right)$ resembling the diagram to the right. 
In particular, in the\\
\begin{picture}(100,5)(5,-5)
\put(360,-100){\epsfig{file=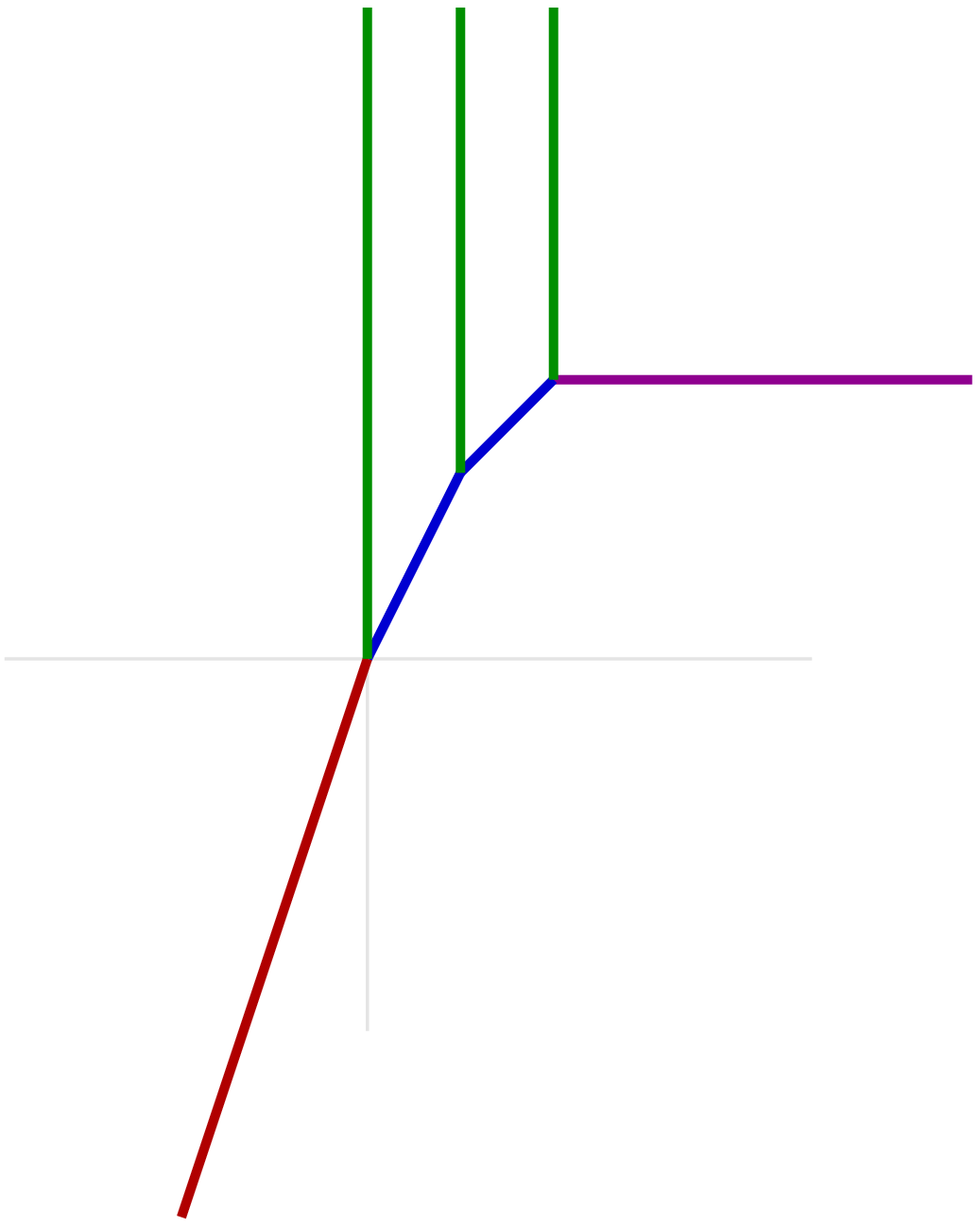,height=2in}} 
\end{picture} 

\vspace{-.4cm} 
\noindent  
notation of Lemma \ref{lemma:vert}, we have\\  
\mbox{}\hspace{2.5cm}$g(x_1)\!:=\!x^3_1-(1+p+p^2)x^2_1+(p+p^2+p^3)x_1-p^3$,\\  
$N\!=\!1$, $t\!=\!4$, and $t'\!=\!3$. \dia  
\end{ex}  

\vspace{-.1cm} 
\noindent 
{\bf Proof of Lemma \ref{lemma:vert}:} By construction, 
$\newt(f)$ lies in a $2$-plane in $\R^{N+1}$\\ 
\scalebox{.99}[1]{and thus, thanks to Kapranov's Theorem, $\ord_p\!\left(
Z^*_{\C_p}(f)\right)$ is the Minkowski}\\ 
\scalebox{.98}[1]{sum of a $1$-dimensional tropical variety and a 
complementary subspace of}\\  
dimension $N-1$. In particular, it suffices to prove the $N\!=\!1$ case. 

The $N\!=\!1$ case follows easily: the ray of type (a) (resp.\ (d)) is 
parallel to the inner normal ray to the edge with vertices 
$(0,1)$ and $(\deg g,0)$ (resp.\ $(0,1)$ and $(0,0)$) of $\newt(f)$. 
The vertical rays correspond to the inner normals corresponding to the lower 
edges of $\newt_p(g)$ (alternatively, the edges of $\newt_p(f)$ not incident 
to $(0,1,0)$). Finally, the ``strips'' are merely the segments connecting 
the points $v$ with $(v,1)$ a lower facet normal of $\newt_p(f)$. \qed 

\medskip 
\noindent 
{\bf Proof of Lemma \ref{lemma:maybetrivial}:} 
Let $n=km+1$ be the number of variables in the system constructed in 
Proposition~\ref{prop:red}. Lemma \ref{lemma:vert} (applied 
to each $y_{i,j}-f_{i,j}$) induces a natural
finite partition of $\R^n$ into half-open slabs of the form
$(-\infty,v_1)\times \R^{n-1}$, $[v_\ell,v_{\ell+1})\times \R^{n-1}$
for $\ell\!\in\!\{1,\ldots,M-1\}$, or $[v_M,+\infty)\times \R^{n-1}$,
with $M\!\leq\!km(t-1)-1$. (Note, in particular, that the boundaries of the
slabs coming from different $f_{i,j}$ can not intersect, thanks to
tropical genericity.) Note also that within the interior of each slab, any
non-empty intersection of $\ord_p\!\left(Z^*_{\C_p}(f_1)\right),\ldots,
\ord_p\!\left(Z^*_{\C_p}(f_n)\right)$ must be a transversal intersection of 
$n$ hyperplanes, thanks to tropical genericity.

In particular, for the left-most (resp.\ right-most) slab,
we obtain a transversal intersection of $n-1$ type (a)
(resp.\ type (d)) $(n-1)$-cells
coming from $f_1,\ldots,f_{n-1}$ (in the notation of Lemma
\ref{lemma:vert}) and an $(n-1)$-cell of $\ord_p\!\left(Z^*_{\C_p}(f_n)
\right)$. Each $(n-1)$-cell of $\ord_p\!\left(Z^*_{\C_p}(f_n)\right)$, by 
definition, is dual to an edge of the lower hull of $\newt_p(f_n)$. So there 
are no more than $\binom{k}{2}$ such $(n-1)$-cells. Thus, there are at most
$\binom{k}{2}$ intersections of
$\ord_p\!\left(Z^*_{\C_p}(f_1)\right),\ldots,\ord_p\!\left(Z^*_{\C_p}(f_n)
\right)$ occuring in the interior of the left-most (resp.\ right-most) slab.

Similarly, the number of intersections of
$\ord_p\!\left(Z^*_{\C_p}(f_1)\right),\ldots,\ord_p\!\left(Z^*_{\C_p}(f_n)
\right)$ occuring in any other slab interior is $\binom{k}{2}$. Also, within 
any of the $M+1$ slab boundaries, there are clearly at most $\binom{k}{2}$
intersections of $\ord_p\!\left(Z^*_{\C_p}(f_1)\right),\ldots,
\ord_p\!\left(Z^*_{\C_p}(f_n)\right)$.

So we obtain no more than $\binom{k}{2}(2+M+M+1)
\!\leq\!\binom{k}{2}(2km(t-1)+1)\!=\!O(k^3mt)$ intersections for 
the underlying $p$-adic tropical varieties and we are done. \qed

\section{Proving Theorem \ref{thm:koi}}  
\label{sec:koi} 
The proof is a fairly straightforward application of a result 
from \cite{koiran}, which we paraphrase in Theorem \ref{thm:hit} 
below. However, let us first review some background. 

Recall that the {\em counting hierarchy} $\ch$ is a hierarchy of
complexity classes built on top of the counting class $\sharpp$;
it contains the entire polynomial hierarchy $\ph$ and 
is contained in $\pspa$. A detailed understanding of $\ch$ is not necessary 
here since we will need only one fact (Theorem \ref{thm:hit} below) related 
to $\ch$. The curious reader can consult
\cite{burgtau,koiran} and the references therein for more information on the
counting hierarchy.

\begin{dfn} 
A {\em hitting set} $H$ for a family $\mathcal{F}$ of polynomials
is a finite set of points such that, for any $f\!\in\!\mathcal F\setminus
\{0\}$, there is at least one $x\!\in\!H$ such that $f(x)\!\neq\!0$.
Also, a {\em $\ch$-algebraic number generator} is 
a sequence of polynomials $G\!:=\!(g_i)_{i\in\N}$ satisfying the following 
conditions: 
\begin{enumerate} 
\item{There is a positive integer $c$ such that 
we can write $g_i(x_1)\!:=\!\sum^{i^c}_{\alpha=0} a(\alpha,i)x^{\alpha}_1$, 
with $a(\alpha,i)\!\in\!\Z$ of absolute value no greater than $2^{i^c}$,  
for all $i$.}  
\item{The language $L(G)\!:=\!\left\{(\alpha,i,j,b)\; | \;  \text{ the } 
j^{\text{\underline{th}}} \text{ bit of } a(\alpha,i) \text{ is equal to } 
b\right\}$ is in $\ch$. \dia }  
\end{enumerate} 
\end{dfn} 

\noindent 
Hitting sets are sometimes called {\em correct test sequences}, as 
in \cite{heintzschnorr}. In particular, the deterministic construction of 
hitting sets is equivalent to the older problem of deterministic 
identity testing for polynomials given in the {\em black-box} model. 

The main technical fact we'll need now is the following:
\begin{thm} (See \cite[Thm.\ 7]{koiran}.) \label{thm:hit}
Let $G\!:=\!(g_{i})$ be a $\ch$-algebraic number generator and let $Z(G,m)$
be the set of all roots of the polynomials $g_{i}$ for all $i \leq m$. 
If there is a polynomial $p$ such that $Z(G,p(kmt))$ is a hitting set for
$\sps(k,m,t)$ then the permanent of 
\scalebox{.91}[1]{$n\times n$ matrices cannot be computed by
constant-free, division-free arithmetic circuits of size $n^{O(1)}$. \qed}  
\end{thm}

\noindent 
The last result shows that the construction of explicit hitting sets
of polynomial size for sums of products of sparse polynomials implies a lower 
bound for the permanent. Note that the conclusion of the theorem holds under a 
somewhat weaker hypothesis (see~\cite{koiran} for details).

\medskip 
\noindent 
{\bf Proof of Theorem~\ref{thm:koi}:} 
By assumption, there is a constant $c\!\geq\!1$ such that any 
$f\!\in\!\sps(k,m,t)$
has at most $(1+kmt)^c$ integer roots that are powers of $p$. (The 
hypothesis of Theorem \ref{thm:koi} is thus in fact stronger than the 
preceding statement.) The set $S_f$ therefore forms a polynomial-size hitting 
set for $f$. By Theorem~\ref{thm:hit}, it just remains to check that the 
sequence of polynomials $(x_1-p^i)_{i \in \N}$ forms a $\ch$-algebraic number 
generator. We must therefore show that the following problem belongs to $\ch$: 
given two integers $i$ and $j$ in binary notation,
compute the $j$-th bit of $p^i$. 
Note that this problem would be solvable in polynomial time if 
$i$ was given in unary notation (by performing the $i-1$ multiplications in 
the most naive way). To deal with the binary notation underlying our 
setting, we apply Theorem 3.10 of~\cite{burgtau}: iterated multiplication of 
exponentially many integers can be done within the counting hierarchy.
Here we have to multiply together exponentially many (in the binary size of 
$i$) copies of the same integer $p$. 
We note that Theorem 3.10 of~\cite{burgtau} applies to a very wide class
of integer sequences: the numbers to be multiplied must  be computable
in the counting hierarchy. In our case we only have to deal with a 
constant sequence (consisting of $i$ copies of $p$) so the elements of
this sequence are computable in polynomial time (and actually in constant time
since $p$ is constant). So we are done. \qed 

\begin{rem} 
\label{rem:stronger} 
From our proof we also obtain that the set $S_f$ from the 
statement of Theorem \ref{thm:koi} can be replaced by 
$S'_f:=\{e\; | \; f\!\left(p^e \right)\!=\!0\}$. 
Since we clearly have $S'_f\!\subseteq\!S_f$, we thus obtain a strengthening 
of Theorem \ref{thm:koi}. \dia 
\end{rem} 

It is interesting to note that even a weakly exponential upper bound 
on the number of valuations would still suffice to prove new hardness results 
for the permanent: from the development of 
Sections 5 and 6 of \cite{koiran}, and our development here, 
one has the following fall-back version of Theorem \ref{thm:koi}. 
\begin{thm} 
\label{thm:closer} 
Suppose that there is a prime $p$ with the following property:
For all $k,m,t\!\in\!\N$ and $f\!\in\!\sps(k,m,t)$,
we have that the cardinality of\\
\mbox{}\hfill $\displaystyle{S'_f:=\left\{e\!\in\!\N \; | \; 
f\!\left(p^e\right)\!=\!0 \right\}}$ 
\hfill\mbox{}\\
is $2^{(kmt)^{o(1)}}$. Then the permanent of $n\times n$ matrices cannot be
computed by polynomial size depth $4$ circuits using polynomial size 
integer constants. \qed   
\end{thm} 

\noindent 
While the conclusion is weaker than that of Theorem \ref{thm:koi}, 
the truth of the hypothesis of Theorem \ref{thm:closer} nevertheless 
yields a hitherto unknown complexity lower bound for the permanent. 
\begin{rem} 
\label{rem:best} 
One can in turn weaken the hypothesis of Theorem \ref{thm:closer} even 
further --- by allowing dependence on the coefficients and degrees of the 
underlying $f_{i,j}$ in the definition of $\sps(k,m,t)$ ---  and still obtain 
the same conclusion. This can be formalized via \cite[Dfn.\ 2.6]{adelic} 
and the development of \cite[Sec.\ 3]{koiran}. \dia 
\end{rem}   

\section{Proving Theorem \ref{thm:upper}}  
\label{sec:upper} 
For the sake of disambiguation, let us first recall the
following basic definition from linear algebra. 
\begin{dfn}
Fix any field $K$. We say that a matrix $E\!=\![E_{i,j}]\!\in\!K^{m\times n}$
is in {\em reduced row echelon form} if and only if the following conditions
hold:
\begin{enumerate}
\item{The left-most nonzero entry of each row of $E$ is $1$, called the
{\em leading $1$} of the row. }
\item{Every leading $1$ is the unique nonzero element of its column.}
\item{\scalebox{.95}[1]{The index $j$ such that $E_{i,j}$ is a leading $1$ 
of row $i$ is a strictly increasing function of $i$. \dia} }
\end{enumerate}
\end{dfn} 

\noindent 
Note that in Condition (1), we allow a row to consist entirely of 
zeroes. Also, by Condition (3), all rows below a row of zeroes must 
also consist solely of zeroes. For example, the matrix 
\scalebox{.5}[.5]{$\begin{bmatrix} \text{\fbox{$1$}} & 0 & 0 & 3  & 0 &  2 \\ 
                 0 & 0 & \text{\fbox{$1$}} & 12 & 0 & -5 \\ 
                 0 & 0 & 0 & 0  & \text{\fbox{$1$}} &  7 \\ 
                 0 & 0 & 0 & 0  & 0 &  0 \end{bmatrix}$} 
is in reduced row echelon form, and we have boxed the leading $1$s. 

By {\em Gauss-Jordan Elimination} we mean the well-known classical algorithm
that, given any matrix $M\!\in\!K^{m\times n}$, yields the factorization
$UM\!=\!E$ with $U\!\in\!\glm(K)$ and $E$ in reduced row echelon form 
(see, e.g., \cite{prasolov,strang}).
In what follows, we use $(\cdot)^\top$ to denote the operation of 
matrix tranpose. 
\begin{dfn} 
\label{dfn:gauss} 
Given any Laurent polynomials 
$f_1,\ldots,f_r\!\in\!K\!\left[x^{\pm 1}_1,\ldots,x^{\pm 1}_n\right]$ 
with supports contained in a set $A\!=\!\{a_1,\ldots,a_t\}\!\subset\!\Zn$ of 
cardinality $t$,  {\em applying Gauss-Jordan Elimination to} 
$(f_1,\ldots,f_r)$ means the following: (a) we identify the row vector 
$(f_1,\ldots,f_r)$ with the vector-matrix product $(x^{a_1},\ldots,x^{a_t})
C$ where $C\!\in\!K^{t\times r}$ and the entries of $C$ are suitably 
chosen coefficients of the $f_i$, and (b) we replace  
$(f_1,\ldots,f_r)$ by $(g_1,\ldots,g_r)$ where $(g_1,\ldots,g_r)\!=\!(x^{a_1},
\ldots,x^{a_t})E$ and $E^\top$ is the reduced row echelon form of $C^\top$. 
\dia 
\end{dfn} 

\noindent 
Note in particular that the ideals $\langle f_1,\ldots,f_r\rangle$ 
and $\langle g_1,\ldots,g_r\rangle$ are identical. As a concrete\linebreak  
example, one can observe that applying Gauss-Jordan Elimination 
to the pair\linebreak $(x^3-y-1,x^3-2y+2)$ means that one instead works with 
the pair $(x^3-4,-y+3)$. 

We now proceed with the proof of Theorem \ref{thm:upper}. 
In what follows, we set $A\!:=\!\bigcup_i A_i$, $t\!:=\!\#A$, and 
let $F\!:=\!(f_1,\ldots,f_n)$ be any polynomial system with 
$f_i\!\in\!\C_p\!\left[x^{\pm 1}_1,\ldots,x^{\pm 1}_n\right]$ and 
$\supp(f_i)\!\subseteq\!A_i$ for all $i$. 

\subsection{Proving Assertions (0) and (1)} 
Assume $t\!\leq\!n+1$. If any $f_i$ is a single monomial term then 
$Z^*_{\C_p}(F)$ is empty. Also, if any $f_i$ is identically $0$ then 
$\#Z^*_{\C_p}(F)$ is infinite, so (by Proposition \ref{prop:fin}) $\ord_p
\!\left(\#Z^*_{\C_p}(F)\right)\!=\!+\infty$. So we may assume that 
no $f_i$ is identically zero or a monomial term. Also, dividing all the $f_i$ 
by a suitable monomial term, we may assume that $\bO\!\in\!A$. 

Assertion (0) then follows immediately from Proposition \ref{prop:over}.   
So we may now assume that $A$ does {\em not} lie in any $(n-1)$-flat 
(and $t\!=\!n+1$ in particular).  

Our remaining case is then folkloric: by 
Gauss-Jordan Elimination (as in Definition \ref{dfn:gauss}, ordering so that 
the last monomial is $x^{\bO}$), 
we can reduce to the case where each polynomial has $2$ or fewer terms, 
and $\supp(f_i)\cap\supp(f_j)\!=\!\bO$ for all $i\!\neq\!j$. 
In particular, should Gauss-Jordan Elimination not yield the 
preceding form, then some $f_i$ is either identically zero or a monomial term, 
thus falling into one of our earlier cases. 
So assume $F$ is a binomial system with $\supp(f_i)\cap\supp(f_j)\!=\!\bO$ for 
all $i\!\neq\!j$. Since no $n+1$ points of $A$ lie on a hyperplane, 
$\newt(f_1),\ldots,\newt(f_n)$ define $n$ linearly independent 
vectors in $\Rn$. The underlying tropical varieties are then hyperplanes 
intersecting transversally, and the number of valuations is thus clearly $1$. 

So the upper bound from Assertion (1) is proved. 
The final equality follows immediately from the polynomial system 
$(x_1-1,\ldots,x_n-1)$. \qed 

\begin{rem} 
Note that in our proof, Gauss-Jordan Elimination allowed us to replace any 
tropically non-generic $F$ by a new, {\em tropically generic} system with the 
same roots over $\C_p$. Recalling standard height bounds for 
linear equations (see, e.g., \cite{storjophd}), another consequence of our 
proof is that, when $t\!\leq\!n+1$, we can 
decide whether $\#\ord_p\!\left(Z^*_{\C_p}(F)\right)$ is $0$, $1$, or 
$\infty$ in polynomial-time. \dia 
\end{rem} 

\subsection{Proving Assertion (2)} 
\scalebox{.9}[1]{Let us first see an example illustrating a trick underlying 
our proof.}   
\begin{ex} 
\label{ex:central} Consider, for any prime $p\!\neq\!2$, the polynomial 
system 

\vspace{-.5cm} 
\noindent 
\begin{minipage}[b]{.75\linewidth}
\vspace{0pt}
$F:=(f_1,f_2):=\left\{\begin{matrix} 
\mbox{} \ \ \ \ \ \ \ \ \ \ \ px^{21}_2 \ \ \ \ \ \ \ \ -px^{32}_1 
\ \ \ \ \ +p \ \ \ \ \ \ \ \ \ \ +x^9_1 x^{10}_2\\ 
-(p+p^2)x^{21}_2+(p+p^3)x^{32}_1+ p+p^4 + (1+p)x^9_1 x^{10}_2 
\end{matrix}\right.$.  
The tropical varieties $\trop_p(f_1)$ and 
$\trop_p(f_2)$ turn out to be\linebreak identical and 
equal to a polyhedral complex with exactly 
$3$\linebreak $0$-dimensional cells and $6$ $1$-dimensional cells (a 
truncation of which is shown on the right). 
\end{minipage} \hspace{.2cm}
\begin{minipage}[b]{.2\linewidth}
\vspace{0pt}
\epsfig{file=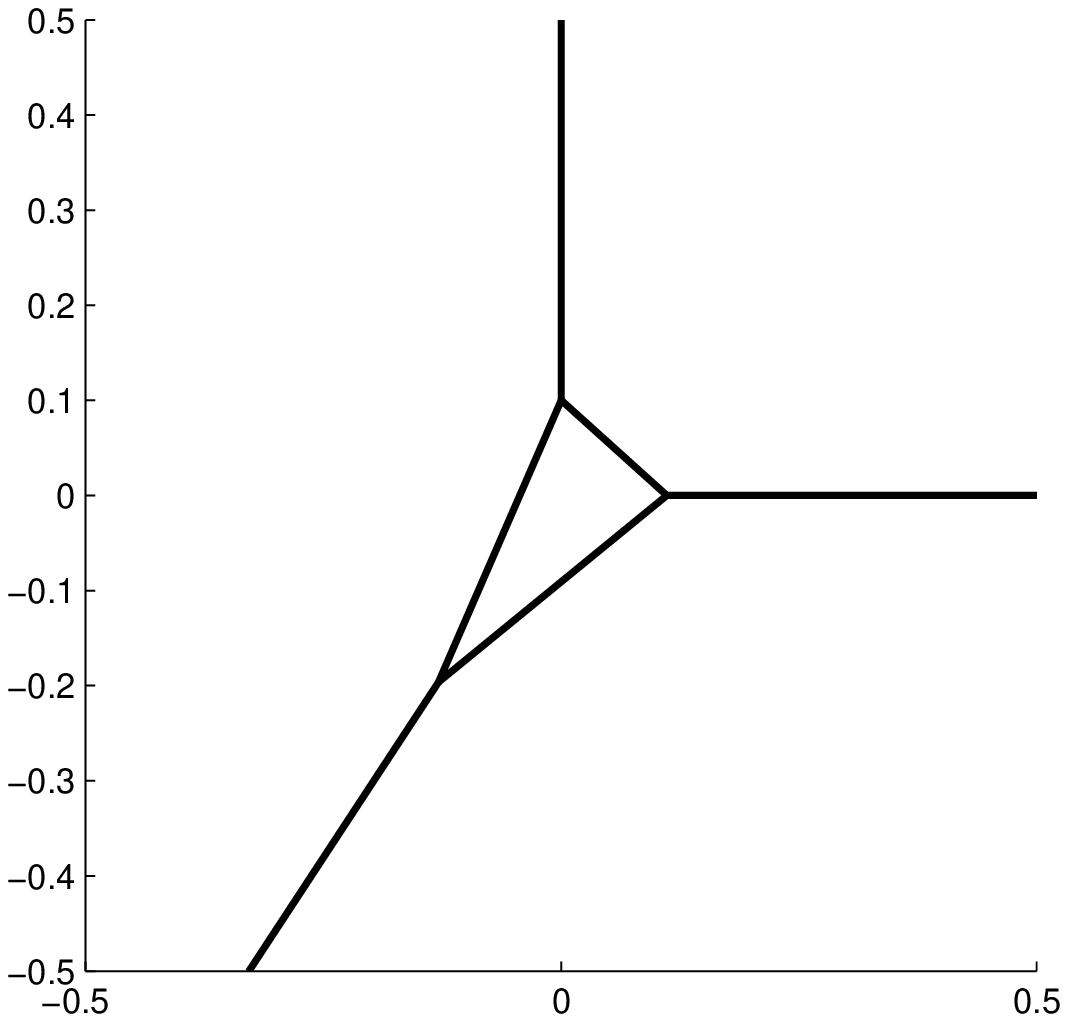,height=1.3in,clip=}
\end{minipage} 

How can we prove that $\ord_p\!\left(Z^*_{\C_p}(F)\right)$ in fact 
has small cardinality? 

\noindent 
\begin{minipage}[b]{.75\linewidth}
\vspace{0pt}
\hspace{.4cm}While it is not hard to apply Bernstein's Theorem 
(as in \cite{bernie}) to see that $F$ has only finitely 
many roots in $(\Cs_p)^2$, there is a simpler approach to 
proving $\ord_p\!\left(Z^*_{\C_p}(F)\right)$ is finite:   
First note that via Gauss-Jordan Elimination (and a suitable 
ordering of monomials), $F$ has the same roots in $(\Cs_p)^2$ as\\ 
\mbox{}\hspace{.1cm}$F^{(1,2)}:=\left(f^{(1,2)}_1,f^{(1,2)}_2\right):=
\left\{\begin{matrix} 
\mbox{} x^{21}_2 \ \ \ \ \ \ \ \ \ \ \ \ \ \ \ +\frac{2+p^2+p^3}{p(p-1)}  
+\frac{2+p+p^2}{p^2(p-1)}x^9_1 x^{10}_2\\ 
\ \ \ \ \ (p+p^3)x^{32}_1 + \frac{2+p+p^3}{p(p-1)} + \frac{2(1+p)}{p^2(p-1)}
x^9_1 x^{10}_2 
\end{matrix}\right.$. 
\end{minipage} 
\begin{minipage}[b]{.2\linewidth}
\vspace{0pt}
\raisebox{-.3cm}{\epsfig{file=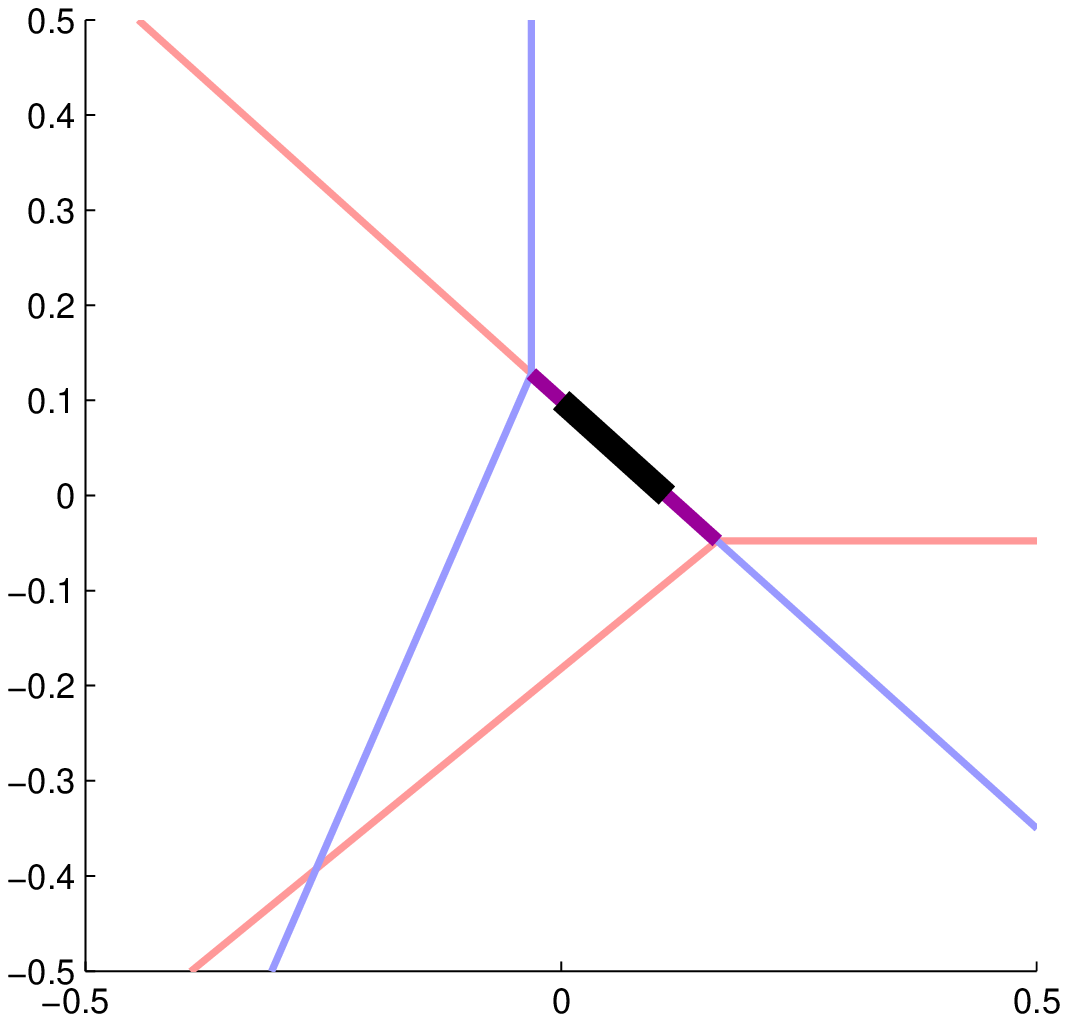,height=1.5in,clip=}}
\end{minipage}\\ 
We then obtain that the tropical varieties 
$\trop_p\!\left(f^{(1,2)}_1\right)$ and $\trop_p\!\left(f^{(1,2)}_2\right)$ 
intersect (in a small interval) along a {\em single} $1$-dimensional 
cell, drawn more thickly, as shown to the right of the definition of 
$F^{(1,2)}$. (The intersection $\trop(f_1)\cap\trop(f_2)\cap
\trop\!\left(f^{(1,2)}_1\right)\cap
\trop\!\left(f^{(1,2)}_2\right)$ is drawn still more thickly.) From the 
definition of $\trop_p(\cdot)$, it is not hard to check that the degenerately 
intersecting $1$-cells of $\trop_p\!\left(f^{(1,2)}_1\right)$ and
$\trop_p\!\left(f^{(1,2)}_2\right)$ correspond to parallel lower edges of 
$\newt_p\!\left(f^{(1,2)}_1\right)$ and $\newt_p\!\left(f^{(1,2)}_2\right)$, 
which in turn correspond to the binomials 
$\frac{2+p^2+p^3}{p(p-1)}+\frac{2+p+p^2}{p^2(p-1)}x^9_1 x^{10}_2$ 
and $\frac{2+p+p^3}{p(p-1)} + \frac{2(1+p)}{p^2(p-1)}
x^9_1 x^{10}_2$. (Note that the intersecting $1$-cells of the 
$\trop\!\left(f^{(1,2)}_i\right)$ are each perpendicular to the resulting 
Newton polytopes of the preceding binomials.) 

\noindent 
\begin{minipage}[b]{.8\linewidth}
\vspace{0pt}
\hspace{.4cm}So to contend with this remaining degenerate intersection, we 
simply apply Gauss-Jordan Elimination with the monomials ordered so that the 
aforementioned pair of binomials becomes a pair of monomials. More precisely, 
we obtain that $F$ has the same roots in $(\Cs_p)^2$ as\\ 
\mbox{}\hfill $F^{(3,4)}:=\left(f^{(3,4)}_1,f^{(3,4)}_2\right):=
\left\{\begin{matrix} 
-\frac{2}{p(p-1)}x^{21}_2 \ + \frac{2+p+p^2}{p(p^2-1)} x^{32}_1 
 +1 \ \ \ \ \ \ \ \ \ \ \ \ \mbox{}\\  
\mbox{} \frac{2-p+p^2}{p(p-1)}x^{21}_2 + \frac{2+p^2+p^3}{p^2-1} x^{32}_1 
\ \ \ \ \ +x^9_1 x^{10}_2 \end{matrix} \right. .$\hfill\mbox{}\\  
\end{minipage} 
\begin{minipage}[b]{.12\linewidth}
\vspace{0pt}
\epsfig{file=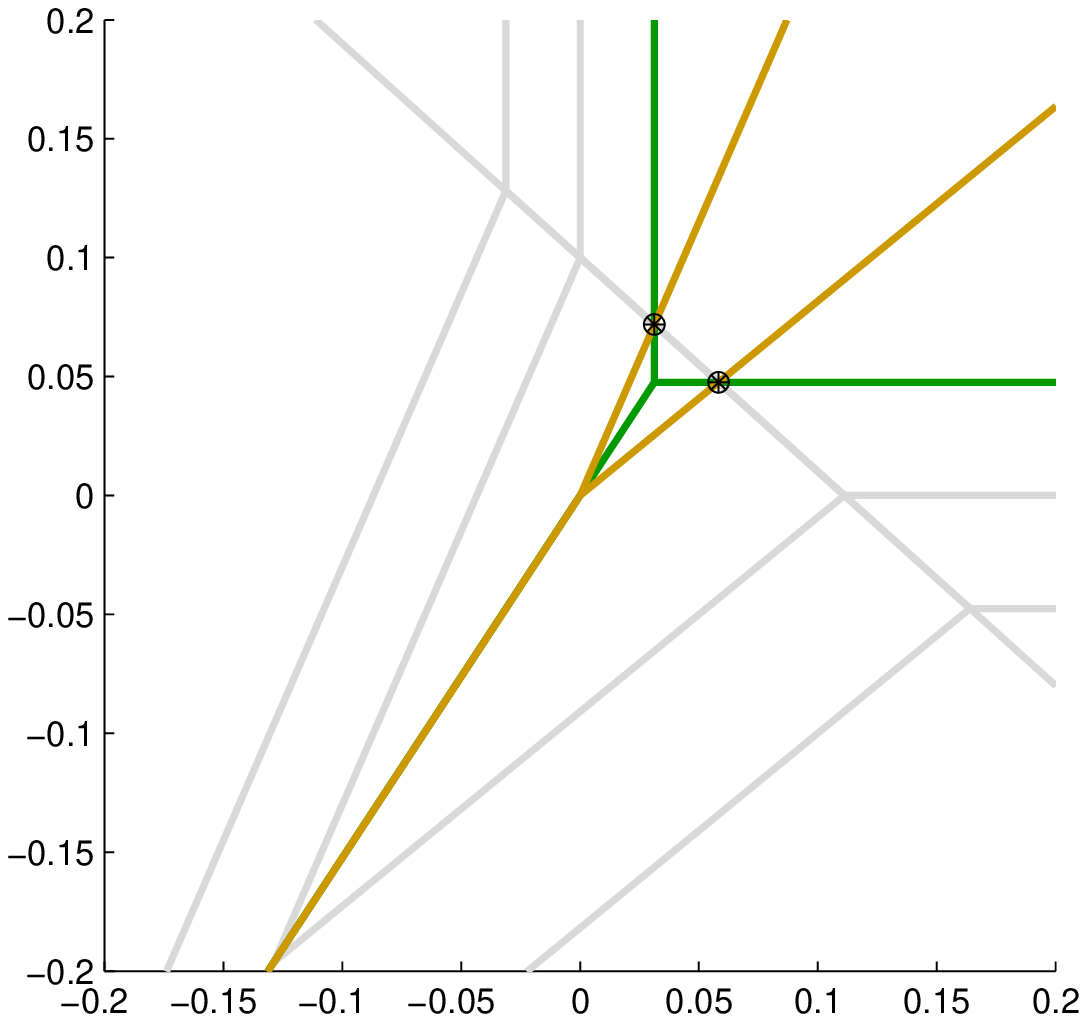,height=1.25in,clip=} 
\end{minipage} 

\noindent 
$\trop_p\!\left(f^{(3,4)}_1\right)$ and $\trop_p\!\left(f^{(3,4)}_2\right)$ 
then intersect transversally precisely within the overlapping $1$-cells of
$\trop_p(f_1)\cap\trop_p(f_2)$ and $\trop_p\!\left(f^{(1,2)}_1\right)\cap
\trop_p\!\left(f^{(1,2)}_2\right)$. In particular, the intersection of the 
tropical varieties of the $f_i$, $f^{(1,2)}_i$, and $f^{(3,4)}_i$ consists of 
exactly $2$ points: $\left(\frac{1}{32},\frac{23}{320}\right)$ and 
$\left(\frac{11}{189},\frac{1}{21}\right)$. So, thanks to Kapranov's 
Theorem, $\ord_p\!\left(Z^*_{\C_p}(F)\right)$ in fact has cardinality at most 
$2$. \dia  
\end{ex} 

A simple observation used in our example above is the following 
consequence of the basic ideal/variety correspondence. 
\begin{prop} 
\label{prop:int} 
Given any $f_1,\ldots,f_r\!\in\!\C_p\!\left[x^{\pm 1}_1,\ldots,x^{\pm 1}_n
\right]$ and $F\!:=\!(f_1,\ldots,f_r)$, we have 
$\ord_p\left(Z^*_{\C_p}(F)\right)\!\subseteq\!\bigcap\limits_{f\in 
\langle f_1,\ldots,f_r\rangle} \ord_p\!\left(Z^*_{\C_p}(f)\right)$, 
where $\langle f_1,\ldots,f_r\rangle\!\subseteq\!\C_p\!\left[x^{\pm 1}_1,
\ldots,x^{\pm 1}_n\right]$ denotes the ideal generated by $f_1,\ldots,f_r$. 
\qed  
\end{prop} 

\noindent 
The reverse inclusion also holds, but is far less trivial to prove  
(see, e.g., \cite{macsturmf}). In particular, Example \ref{ex:central} shows 
us that restricting the intersection to a {\em finite} set of generators 
can sometimes result in strict containment. 

\medskip 
\noindent 
{\bf Proof of Assertion (2) of Theorem \ref{thm:upper}:} 
The case $n\!=\!1$ is immediate from Lemma \ref{lemma:newt} and the example 
$(A_1,f)\!=\!(\{0,1,2\},(x_1-1)(x_1-2))$. So let us assume henceforth that 
$n\!\geq\!2$.  

Recall from last section that $t\!=\!n+2$ 
and $A\!:=\!\bigcup_{\ell}A_\ell$. For any distinct 
$i,j\!\in\!\{1,\ldots,n+2\}$ let us then define 
$F^{(i,j)}$ by applying Gauss-Jordan Elimination,  
as in Definition \ref{dfn:gauss}, where we order monomials   
so that the last exponents are $a_i$ and $a_j$. 
In particular, $F$ and $F^{(i,j)}$ clearly have the same roots in $(\Cs_p)^n$ 
for all distinct $i,j$. 
We will show that {\em every} $F^{(i,j)}$ can be assumed 
to be a trinomial system of a particular form. 

First note that we may assume that 
$\supp\!\left(f^{(n+1,n+2)}_\ell\right)\!\subseteq\!\{a_{r_\ell},a_{s_\ell},
a_{n+2}\}$ for all $\ell$, where $(r_\ell)_\ell$ is a strictly increasing 
sequence of integers in $\{1,\ldots,n\}$ satisfying 
$s_\ell\!\geq\!r_\ell\!\geq\!\ell$ for all $\ell$. This is because, similar to 
our last proof, we may assume that each $f^{(n+1,n+2)}_\ell$ has at least $2$ 
monomial terms, thus implying that $r_n\!\in\!\{n,n+1\}$. In particular, no 
$f^{(n+1,n+2)}_\ell$ can have $4$ or more terms, by the positioning of the 
leading $1$s in reduced row echelon form.  

By dividing by a suitable monomial term, we may assume that $a_{n+2}\!=\!\bO$. 
Also, by Lemma \ref{lemma:mono}, we may assume that 
$a_\ell\!\in\!\Z^{\ell}\times \{0\}^{n-\ell}$ for all $\ell\!\in\!\{1,\ldots,
n\}$. (Our general position assumption on $A$ also implies that the 
$\ell\thth$ coordinate of $a_\ell$ is nonzero.) 
Now, should $f^{(n+1,n+2)}_n$ have exactly $2$ monomial terms, then 
$f^{(n+1,n+2)}_n$ must be of one of the following forms: 
(a) $x^{a_n}+\alpha_n x^{a_{n+1}}$,
(b) $x^{a_n}+\alpha_n x^{a_{n+2}}$, or (c) $x^{a_{n+1}}+\alpha_n x^{a_{n+2}}$,
for some $\alpha_n\!\in\!\Cs_p$. In Case (c), we could then replace
all occurences of $x^{a_{n+1}}$ in $f^{(n+1,n+2)}_1,\ldots,
f^{(n+1,n+2)}_{n-1}$ by a nonzero multiple of $x^{a_{n+2}}$. We would thus 
reduce to the setting of Assertion (1), in 
which case, the maximal finite number of valuations would be $1$. In Cases 
(a) and (b), we obtain either that some $x_i$ vanishes, or that 
$\ord_p x_n$ is a linear function of $\ord_p x_1,\ldots,
\ord_p x_{n-1}$. So we could then reduce to a case one dimension lower. 

So we may assume that 
$\supp\!\left(f^{(n+1,n+2)}_n\right)\!=\!\{a_n,a_{n+1},a_{n+2}\}$, which in 
turn forces $a_\ell\!\in\!\supp\!\left(f^{(n+1,n+2)}_\ell\right)$ 
$\subseteq\!\{a_\ell,a_{n+1},a_{n+2}\}$ for all $\ell\!\in\!\{1,\ldots,n-1\}$. 
Moreover, by repeating the arguments of Cases (a) and (b) above, we may in fact 
assume $\supp\!\left(f^{(n+1,n+2)}_\ell\right)\!=\!\{a_\ell,a_{n+1},a_{n+2}\}$ 
for all $\ell\!\in\!\{1,\ldots,n\}$. 

Permuting indices, we can then repeat the last $3$ paragraphs  
and assume further that, for any distinct $i,j\!\in\!\{1,\ldots,n+2\}$, 
we have 

\medskip 
\noindent 
($\star$) \hspace{.5cm} 
$\supp\!\left(f^{(i,j)}_\ell\right)\!=\!\{a_{k_\ell},a_i,a_j\}$ for 
all $\ell\!\in\!\{1,\ldots,n\}$, where $\{k_\ell\}_\ell\!=\!A\!\setminus
\!\{i,j\}$. 

\medskip 
\noindent 
Let us now fix $(i,j)$ and set $G\!=\!(g_1,\ldots,g_n)\!:=\!F^{(i,j)}$. 
Thanks to ($\star$) and Lemma \ref{lemma:vert} (mimicking the proof of 
Lemma \ref{lemma:maybetrivial}), the $\trop_p(g_i)$ each contain 
a half-plane parallel to a common hyperplane. We then obtain a 
finite partition of $\R^n$ into half-open slabs of a
form linearly isomorphic (over $\Q$) to
$(-\infty,u_1)\times \R^{n-1}$, $[u_\ell,u_{\ell+1})\times \R^{n-1}$
for $\ell\!\in\!\{1,\ldots,m_{i,j}-1\}$, or $[u_{m_{i,j}},+\infty)\times 
\R^{n-1}$, with $m_{i,j}\!\leq\!n$. More precisely, the boundaries of the 
cells of our partition are hyperplanes of the form 
$H^{(i,j)}_\ell\!:=\!\left\{v\!\in\!\Rn\; \left| 
\; (a_i-a_j)\cdot v\!=\!\ord_p\!
\left( \gamma^{(i,j)}_\ell\right)\right. \right\}$ for 
$\ell\!\in\!\{1,\ldots,m_{i,j}\}$, where 
$\gamma^{(i,j)}_\ell$ is a ratio 
of coefficients of $f^{(i,j)}_\ell$. 

In particular, Lemma \ref{lemma:vert} and our genericity hypothesis 
tell us that, within any slab, any non-empty intersection of 
$\ord_p\!\left(Z^*_{\C_p}(f_1)\right),\ldots,
\ord_p\!\left(Z^*_{\C_p}(f_n)\right)$ must be a {\em transversal}  
intersection of $n$ hyperplanes, unless it includes the intersection 
of two or more $H^{(i,j)}_{k}$. So if $G$ is tropically generic,  
we have by Proposition \ref{prop:int} that $\#\ord_p\!\left(Z^*_{\C_p}(F)
\right)\!\leq\!n+1$.

Otherwise, any non-transversal intersection must 
occur within an intersection of slab boundaries $H^{(i,j)}_k$. 
So to finish this case, 
consider $n-1$ more distinct pairs $(i_2,j_2),\ldots,
(i_n,j_n)$, i.e., $i_\ell\!\neq\!j_\ell$ for all $\ell$ and 
$\#\{i_\ell,j_\ell,i_{\ell'},j_{\ell'}\}\!\leq\!3$ for all 
$\ell\!\neq\!\ell'$. 

Just as for $G$, the genericity of the 
exponent set $A$ implies that any non-transversal intersection 
for $\trop_p\!\left(f^{(i_\ell,j_\ell)}_1\right),
\ldots,\trop_p\!\left(f^{(i_\ell,j_\ell)}_n\right)$ must 
occur within the intersection of at least two coincident slab boundaries 
$H^{(i_\ell,j_\ell)}_k$. In particular, we may assume that none of 
$F^{(i_2,j_2)},\ldots,F^{(i_n,j_n)}$ are tropically generic (for 
$\#\ord_p\!\left(Z^*_{\C_p}(F)\right)\!\leq\!n+1$ otherwise). 

By our assumption on the genericity of the exponent set 
$A$, we have that $H^{(i,j)}_{\ell_1},H^{(i_2,j_2)}_{\ell_2},
\ldots,$ 
$H^{(i_n,j_n)}_{\ell_n}$ intersect transversally, 
for any choice of $n$-tuples $(\ell_1,\ldots,\ell_n)$. 
In particular, we have embedded the non-transversal 
intersections of $\ord_p\!\left(Z^*_{\C_p}(F)\right)$ into a 
(finite) intersection of $m_{i,j}\prod^n_{\ell=2} m_{i_\ell,j_\ell}$ many 
tropical varieties. In particular, to count the non-transversal 
intersections, we may assume 
$m_{i,j},m_{i_2,j_2},\ldots,m_{i_n,j_n}\!\leq\!\floor{\frac{n}{2}}$.  

From our earlier observations on slab decomposition, 
there can be at most $n$ intersections occuring away from an 
intersection of slab boundaries (since we are assuming 
$G$ and the $F^{(i_\ell,j_\ell)}$ all fail to be tropically generic). 
The number of distinct points of $\ord_p\!\left(Z^*_{\C_p}(F)\right)$ lying  
in intersections of the form $H^{(i,j)}_{\ell_1}\cap H^{(i_2,j_2)}_{\ell_2}
\cap \cdots \cap H^{(i_n,j_n)}_{\ell_n}$ is no greater than 
$\floor{\frac{n}{2}}^n$. So our upper bound is proved.  

That $\cV_p\!\left(\{\bO,2e_1,e_1+e_2\},\{\bO,2e_1,e_2+e_3\},\ldots,
\{\bO,2e_1,e_{n-1}+e_n\},\{\bO,2e_1,e_n\}\right)\!\geq\!n+1$ 
follows directly from \cite[Thm.\ 1.6]{adelic}. To be more precise, the 
polynomial system\\ 
\scalebox{.84}[1]{ 
$\displaystyle{\left(x_1x_2-
p\left(1+\frac{x^2_1}{p}\right),
x_2x_3-\left(1+p x^2_1\right),
x_3x_4-\left(1+p^3x^2_1\right),
\ldots,x_{n-1}x_n-\left(1+p^{2n-5}x^2_1\right),
x_n-\left(1+p^{2n-3}x^2_1\right)\right)}$} 
has exactly $n+1$ valuation vectors for its roots over $\C_p$, and 
tropical genericity follows directly from \cite[Lemma 3.7]{adelic}. 
The reverse inequality then 
follows from our earlier observations on slab decomposition. 
In particular, via our earlier reductions, Assertions (0) and (1) 
easily\linebreak 
\scalebox{.95}[1]{imply that any $F$ with smaller support has no more than $n$ 
valuation vectors for its roots. \qed}  
\begin{rem} 
We are currently unaware of any examples where 
$\blah_p(A_1,\ldots,A_n)$ is larger than $n+1$. In any event, our 
proof reveals various cases where the number of valuation vectors is 
at most $n+1$.  \dia 
\end{rem} 

\section{Proving Theorem \ref{thm:hard}}  
\label{sec:hard} 
By Theorem \ref{thm:koi} and Proposition \ref{prop:red} it is enough to show 
that, for $n\!:=\!km+1$ and $A_1,\ldots,A_n$ the supports of the 
polynomial system $F$ from the proposition, we have\linebreak  
$\cV_p(A_1,\ldots,A_n)\!=\!(kmt)^{O(1)}$. 
By Lemma \ref{lemma:maybetrivial} we are done. \qed 

\section{Proving Theorem \ref{thm:mult}} 
\label{sec:mult} 
First note that we can divide our equations by a suitable monomial 
term so that $\bO\!\in\!A$. 

The case where $A$ has cardinality $n+1$ can be easily 
handled just as in the proof of Theorem \ref{thm:upper}: 
$F$ can be reduced to a binomial system via Gauss-Jordan Elimination,  
and then by Lemma \ref{lemma:mono} we can easily reduce to a {\em triangular} 
binomial system. In particular, all the roots of $F$ in $(\Ks)^n$ are 
non-degenerate and thus have multiplicity $1$, so the sharpness of the bound is 
immediate as well. 

So let us now assume that $A$ has cardinality $n+2$. 
By Gauss-Jordan Elimination and a monomial change of 
variables again, we may assume that $F$ is of the form\\ 
\mbox{}\hfill $\left(x^{a_1} - \alpha_1 - x^{a_{n+1}}/c, \ldots, 
x^{a_n} - \alpha_n - x^{a_{n+1}}/c\right)$\hfill\mbox{}\\  
for some $c\!\in\!\Ks$. 

Consider now the matrix $\hat{A}$ obtained by appending a rows 
of $1$s to the matrix with columns $\bO,a_1,\ldots,a_{n+1}$. By 
construction, $\hat{A}$ has right-kernel generated by a single 
vector $b\!=\!(b_0,\ldots,b_{n+1})\!\in\!\Z^{n+2}$ with 
{\em no} zero coordinates. 
So the identity $1^{b_0}\left(x^{a_1}\right)^{b_1}\cdots\left(
x^{a_{n+1}}\right)^{b_{n+1}}\!=\!1$ clearly holds for any 
$x\!\in\!(\Ks)^n$. Letting $u\!:=\!x^{a_{n+1}}$ we then clearly 
obtain a bijection between the roots of $F$ in $(\Ks)^n$ and the roots of \\ 
\mbox{}\hfill $g(u)\!:=\!u^{b_{n+1}}\left(\prod^n_{i=1}(\alpha_i +u)^{b_i}
\right)-C$ \hfill \mbox{}\\ 
where 
$C\!:=\!c^{b_1+\cdots+b_n}$. 
Furthermore, intersection multiplicity is preserved under this univariate 
reduction since each $x_i$ is a radical of a linear function of a root of 
$g$. We thus need only determine the maximum intersection multiplicity 
of a root of $g$ in $\Ks$. 

Since the multiplicity of a root $\zeta$ over a field of 
characteristic $0$ is characterized by the derivative of 
least order not vanishing at $\zeta$, let us suppose, to 
derive a contradiction, that $f(\zeta)\!=\!f'(\zeta)\!=\cdots=
\!f^{(n+1)}(\zeta)\!=\!0$, i.e., $\zeta$ is a root of multiplicity 
$\geq\!n+2$. An elementary calculation then reveals that we must have 
\begin{eqnarray*} 
& \frac{b'_1}{\alpha'_{1}+\zeta}+\cdots+\frac{b'_{m+1}}{\alpha'_{m+1}+\zeta} 
& = 0\\ 
& \vdots & \\ 
& \frac{b'_1}{(\alpha'_{1}+\zeta)^{n+1}}+\cdots+\frac{b'_{m+1}}
{(\alpha'_{m+1}+\zeta)^{n+1}} & = 0 
\end{eqnarray*} 

\noindent 
where $m\!\leq\!n$, the $\alpha'_i$ are distinct 
and comprise all the $\alpha_i$, $\alpha'_{m+1}\!=\!0$, $b'_i\!:=\!\sum
\limits_{\alpha_{j}=\alpha'_i} b_{j}$, $b'_{m+1}\!:=\!b_{n+1}$, 
we set $\alpha'_{m+1}\!:=\!0$, and $\zeta\!\not\in\!\{-\alpha_i\}$. 
In other words, $[b'_1,\ldots,b'_{m+1}]^\top$ is a right-null vector of a 
Vandermonde matrix with non-vanishing determinant. Since 
$[b'_1,\ldots,b'_{m+1}]$ has nonzero coordinates, we thus obtain a 
contradiction. So our upper bound is proved. 

To prove that our final bound is tight, let $\zeta_1,\ldots,\zeta_{n+1}$ 
denote the (distinct) $(n+1)\stst$ roots of unity in $K$ and 
set $g(u)\!:=\!\displaystyle{u\left(\prod^n_{i=1} (u+\zeta_{n+1}-\zeta_i)
\right)+1}$. Since $g(u-\zeta_{n+1})\!=\!u^{n+1}$, it is clear that $g$ has 
$-\zeta_{n+1}$ as a root of multiplicity $n+1$. Furthermore, $g$ is 
nothing more than the univariate reduction argument of our proof 
applied to the system 
\begin{eqnarray*} 
\theta x_1 & = & \zeta_{n+1}-\zeta_i +\frac{1}{x_1\cdots x_n} \\ 
& \vdots  &  \\ 
\theta x_n & = & \zeta_{n+1}-\zeta_n +\frac{1}{x_1\cdots x_n} 
\end{eqnarray*} 

\noindent 
where $\theta$ is any $n\thth$ root of $-1$. \qed  

\section*{Acknowledgements} 
We thank Henry Cohn for his wonderful hospitality at 
Microsoft Research New England (where Rojas presented a preliminary version 
of Lemma \ref{lemma:maybetrivial} on July 30, 2012), and  
Kiran Kedlaya and Daqing Wan for useful $p$-adic discussions. We also 
thank Mart\'\i{}n Avendano, Bruno Grenet, and Korben Rusek for 
useful discussions on an earlier version of Lemma \ref{lemma:vert}, and 
Jeff Lagarias for insightful comments on an earlier version of this paper. 

Most importantly, however, we would like to congratulate Mike Shub on 
his 70$\thth$ birthday: he has truly blessed us with his friendship and his 
beautiful mathematics. We hope this paper will serve as a small but 
nice gift for Mike. 

\bibliographystyle{amsalpha}

\begin{thebibliography}{A}   


\bibitem[Art67]{artin} Artin, Emil, {\it Algebraic Numbers and Algebraic 
Functions,} Gordon and Breach, New York, 1967. 



\bibitem[Ber75]{bernie} Bernshtein, David N., {\it ``The Number of
Roots of a System of Equations,"} Functional Analysis and its Applications
(translated from Russian), Vol. 9, No. 2, (1975), pp.\ 183--185.

\bibitem[BCSS98]{bcss} Blum, Lenore; Cucker, Felipe; Shub, Mike; and
Smale, Steve, {\it Complexity and Real Computation,} Springer-Verlag, 1998. 

\bibitem[BC76]{bocook} Borodin, Alan and Cook, Steve, {\it ``On the number 
of additions to compute specific polynomials,''} SIAM Journal on 
Computing, 5(1):146--157, 1976.  

\bibitem[Bra39]{brauer} Brauer, Alfred, {\it ``On addition chains,''} 
Bull.\ Amer.\ Math.\ Soc.\ 45, (1939), pp.\ 736--739.  


\bibitem[B\"ur00]{burgcook} B\"urgisser, Peter, {\it ``Cook's 
versus Valiant's Hypothesis,''} Theor.\ Comp.\ Sci., 
235:71--88, 2000.



\bibitem[B\"ur09]{burgtau} \underline{\hspace{\burg}}, {\it ``On defining 
integers and proving arithmetic circuit lower bounds,''} 
Computational Complexity, 18:81--103, 2009. 

\bibitem[BLMW11]{jml} B\"urgisser, Peter; Landsberg, J.\ M.; 
Manivel, Laurent; and Weyman, Jerzy, {\it ``An Overview of Mathematical 
Issues Arising in the Geometric Complexity Theory Approach to 
$\mathbf{VP}\neq\mathbf{VNP}$,''} SIAM J.\ Comput.\ {\bf 40}, pp.\ 1179-1209, 
2011.  

\bibitem[Che04]{cheng} Cheng, Qi, {\it ``Straight Line Programs and
Torsion Points on Elliptic Curves,''} Computational Complexity,
Vol.\ 12, no.\ 3--4 (sept.\ 2004), pp.\ 150--161.



\bibitem[EKL06]{kapranov} Einsiedler, Manfred; Kapranov, Mikhail;
and Lind, Douglas, {\it ``Non-archimedean amoebas and tropical varieties,''}
Journal f\"ur die reine und angewandte Mathematik (Crelles Journal), Vol.
2006, no.\ 601, pp.\ 139--157, December 2006.

\bibitem[Ful08]{ifulton} Fulton, William, {\it
Intersection Theory,} 2$\nd$ ed., Ergebnisse der Mathematik und
ihrer Grenzgebiete 3, {\bf 2}, Springer-Verlag, 2008.

\bibitem[vzGG03]{compualg} von zur Gathen, Joachim and Gerhard, 
J\"urgen, {\it ``Modern Computer Algebra,''} 2nd ed., Cambridge University 
Press, 2003. 


\bibitem[Gou03]{gouvea} Gouv\^ea, Fernando Q., {\it $p$-adic 
Numbers,} Universitext, 2nd ed., Springer-Verlag, 2003. 

\bibitem[HS82]{heintzschnorr} Heintz, Joos and Schnorr, C.-P., 
{\it ``Testing polynomials which are easy to compute,''} in 
Logic and Algorithmic 
(international symposium in honor of Ernst Specker), pp.\ 237--254, 
monograph no.\ 30 of L'Enseignement Math\'ematique, 1982. 

\bibitem[Her56]{hermite} Hermite, Charles, {\it ``Sur le nombres des racines
d'une \'equation alg\'ebrique compris\'e entre des limites donn'es,"} J.\
Reine Angew.\ Math.\ 52 (1856) 39--51; also: Oeuvres, Vol. I
(Gauthier-Villars, Paris, 1905) pp.\ 397--414; English translation: P.\ C.\
Parks, Internat.\ J.\ Control 26 (1977), pp.\ 183--195.




\bibitem[Kat07]{katok} Katok, Svetlana, {\it $p$-adic Analysis 
Compared with Real,} Student Mathematical Library, vol.\ 37, 
American Mathematical Society, 2007.  


\bibitem[Koi11]{koiran} Koiran, Pascal, {\it ``Shallow Circuits with 
High-Powered Inputs,''} in Proceedings of 
Innovations in Computer Science (ICS 2011, Jan.\ 6--9, 2011, Beijing China), 
Tsinghua University Press, Beijing. 




\bibitem[Lip94]{lipton} Lipton, Richard, {\it ``Straight-line complexity and 
integer factorization,''} Algorithmic number theory (Ithaca, NY, 1994), pp.\ 
71--79, Lecture Notes in Comput.\ Sci., 877, Springer, Berlin, 1994. 

\bibitem[MS12]{macsturmf} Maclagan, Diane and Sturmfels, Bernd,
{\it Introduction to Tropical Geometry,} in progress.

\bibitem[Maz78]{mazur} Mazur, Barry, {\it ``Rational Isogenies of Prime 
Degree,''} Invent.\ Math., 44, 1978. 

\bibitem[dMS96]{svaiter} de Melo, W.\ and Svaiter, B.\ F.,
{\it ``The cost of computing integers,''} Proc.\ Amer.\ Math.\ Soc.\
{\bf 124} (1996), pp.\ 1377--1378.

\bibitem[Mer96]{merel} Merel, Loic, {\it ``Bounds for the torsion of 
elliptic curves over number fields,''} Invent.\ Math., 
124(1--3):437--449, 1996. 

\bibitem[Mor97]{moreira} T.\ de Araujo Moreira, Gustavo, 
{\it ``On asymptotic estimates for arithmetic cost functions,''} 
Proccedings of the American Mathematical Society, Vol.\ 125, no.\ 2, 
Feb.\ 1997, pp.\ 347--353.  


\bibitem[Par99]{parent} Parent, Philippe, {\it ``Effective Bounds for the 
torsion of elliptic curves over number fields,''} J.\ Reine Angew.\ Math, 
508:65--116, 1999. 

\bibitem[PR13]{adelic} Phillipson, Kaitlyn and Rojas, J.\ Maurice,
{\it ``Fewnomial Systems with Many Roots, and an Adelic Tau
Conjecture,''} in proceedings of Bellairs workshop on tropical and
non-Archimedean geometry (May 6--13, 2011, Barbados), 
Contemporary Mathematics, vol.\ 605, pp.\ 45--71, AMS Press, to appear. 

\bibitem[Pra04]{prasolov} Prasolov, V.\ V., {\it Problems and Theorems 
in Linear Algebra,} translations of mathematical monographs, vol.\ 134, 
AMS Press, 2004. 


\bibitem[Rob00]{robert} Robert, Alain M.,
{\it A course in $p$-adic analysis,}
Graduate Texts in Mathematics, 198, Springer-Verlag, New York, 2000. 




\bibitem[Sch84]{schikhof} Schikhof, W.\ H., {\it Ultrametric Calculus, An
Introduction to $p$-adic Analysis,} Cambridge Studies in
Adv.\ Math.\ 4, Cambridge Univ. Press, 1984.

\bibitem[Ser79]{serre} Serre, Jean-Pierre, {\it Local fields,} 
Graduate Texts in Mathematics, 67, Springer-Verlag, New York-Berlin, 1979. 

\bibitem[Shu93]{shub} Shub, Mike, {\it ``Some Remarks
on B\'ezout's Theorem and Complexity Theory,''} {}From
Topology to Computation: Proceedings of
the Smalefest (Berkeley, 1990), pp.\ 443--455, Springer-Verlag, 1993. 

\bibitem[Sma98]{21a} Smale, Steve, {\it
``Mathematical Problems for the Next Century,''}
Math.\ Intelligencer  20  (1998),  no.\ 2, pp.\ 7--15.

\bibitem[Sma00]{21b} \underline{\hspace{\sma}}, {\it
``Mathematical Problems for the Next Century,''}
Mathematics: Frontiers and Perspectives, pp.\ 271--294, Amer.\ Math.\ Soc.,
Providence, RI, 2000. 

\bibitem[Smi97]{smirnov} Smirnov, A.\ L., {\it ``Torus Schemes Over a
Discrete Valuation Ring,''}  St.\ Petersburg Math.\ J.\ {\bf 8} (1997),
no.\ 4, pp.\ 651--659.

\bibitem[Sto00]{storjophd} Storjohann, Arne, {\it ``Algorithms for
matrix canonical forms,''} doctoral dissertation, Swiss Federal
Institute of Technology, Zurich, 2000. 

\bibitem[Str09]{strang} Strang, Gilbert, {\it Introduction to Linear 
Algebra,} 4$\thth$ edition, Wellesley-Cambridge Press, 2009. 

\bibitem[Val79]{valiant} Valiant, Leslie G., {\it ``The complexity of 
computing the permanent,''} Theoret.\ Comp.\ Sci., 8:189--201, 1979. 

\bibitem[Wei63]{weiss} Weiss, Edwin, \emph{Algebraic Number Theory}, 
McGraw-Hill, 1963.

\bibitem[Zie95]{ziegler} Ziegler, Gunter M., {\em Lectures on
Polytopes}, Graduate Texts in Mathematics, Springer
Verlag, 1995.

\end{thebibliography}

\end{document}